\documentclass[11pt,twoside]{article}
\usepackage{amsmath}
\usepackage{amssymb}
\usepackage{latexsym}

\usepackage[T1]{fontenc}

\usepackage[all,2cell]{xy}
\UseAllTwocells
\newdir{ >}{{}*!/-5pt/@{>}}

\renewcommand{\thesection}{\arabic{section}}

\usepackage{titlesec}

\titleformat{\section}[block]{\bfseries\filcenter}%
{{\upshape\S\,\thesection\enspace}}{.5em}{}

\usepackage[labelsep=period]{caption}[2006/03/21]

\newtheorem{df}{Definition}[section]
\newtheorem{thm}[df]{Theorem}
\newtheorem{prop}[df]{Proposition}
\newtheorem{cor}[df]{Corollary}
\newtheorem{ex}[df]{Example}

\newtheorem{lem}[df]{Lemma}

\newcommand{\pf}{\noindent{\sc Proof.}\ }

\newcommand{\pfend}%
{\hspace*{\fill}\lower3pt\hbox{$\Box$}\medskip}

\newcommand{\D}{\mathop{\raise0.1ex\hbox{$\mathfrak{D}$}}} 

\newcommand{\chigh}{{\raise1.5pt\hbox{$\chi$}}}
\renewcommand{\phi}{\varphi}
\newcommand{\Ga}{\Gamma}
\newcommand{\Om}{\Omega}

\newcommand{\R}{\mathbb{R}}

\newcommand{\id}{\text{\upshape id}}

\newcommand{\da}{\partial}

\newcommand{\tilzero}{\widetilde 0}         
\newcommand{\tilq}{\widetilde q}
\newcommand{\tilR}{\widetilde R}

\newcommand{\tilalpha}{\widetilde\alpha}
\newcommand{\tilbeta}{\widetilde\beta}

\newcommand{\upa}{\uparrow}

\newcommand{\jd}{\Tilde{J}}
\newcommand{\td}{\Tilde{\Theta}}

\usepackage{mathrsfs}
\newcommand{\mscr}[1]{\mathscr{#1}}

\newcommand{\VBgpd}{$\mscr{VB}$--groupoid}
\newcommand{\VBgpds}{$\mscr{VB}$--groupoids}
\newcommand{\LAvb}{$\mscr{LA}$--vector bundle}
\newcommand{\LAgpd}{$\mscr{LA}$--groupoid}

\newcommand{\LAgpds}{$\mscr{LA}$--groupoids}

\newcommand{\gog}{\mathfrak{g}}

\newcommand{\Goi}{\mathfrak{I}}

\newcommand{\gop}{\mathfrak{p}}
\newcommand{\goR}{\mathfrak{R}}

\newcommand{\cinf}[1]{C^\infty(#1)}

\usepackage{euscript}

\newcommand{\extt}[2]{\mathsf{\Lambda}^{#1}(#2)}

\newcommand{\ld}{\mathfrak{L}}             

\newcommand{\co}{\colon\thinspace}

\renewcommand{\leq}{\leqslant}

\newcommand{\isom}{\cong}

\newcommand{\xclam}{^{\textstyle !}}

\newcommand{\sdp}{\ltimes}
\newcommand{\pds}{\rtimes}

\newcommand{\act}{\mathbin{\hbox{$<\kern-.4em\mapstochar\kern.4em$}}}
\newcommand{\ract}{\mathbin{\hbox{$\mapstochar\kern-.3em>$}}}

\newcommand{\sol}{\bullet}

\newcommand{\Whitney}[1]{\mathbin{\oplus_{#1}}}

\newcommand{\llangle}{\langle\!\langle}
\newcommand{\rrangle}{\rangle\!\rangle}
\newcommand{\lllangle}{\langle\!\langle\!\langle}
\newcommand{\rrrangle}{\rangle\!\rangle\!\rangle}

\usepackage{stmaryrd}

\newcommand{\thbr}[2]%
{\rule[-1pt]{1pt}{10pt}\hspace{2pt} #1,\, #2\hspace{1pt}\rule[-1pt]{1pt}{10pt}}

\newcommand{\newp}[2]%
{\talloblong\hspace{1pt} #1,\, #2\hspace{0.5pt}\talloblong}

\renewcommand{\Bar}[1]{\overline{#1}}
\renewcommand{\Tilde}[1]{\widetilde{#1}}

\newcommand{\minilos}{\medskip%
\centerline{\textbf{\small *\qquad *\qquad *\qquad *\qquad *\qquad}}%
\smallskip}

\newcommand{\gpd}{\rightrightarrows} 

\newcommand{\dan}{\Delta}

\usepackage{wasysym}

\newcommand{\add}[1]{\mathbin{\lower 5pt%
    \hbox{${\stackrel{\textstyle +}{\scriptscriptstyle #1}}$}}}

\newcommand{\less}[1]{\mathbin{\lower 5pt%
    \hbox{${\stackrel{\textstyle -}{\scriptscriptstyle #1}}$}}}

\newcommand{\by}[1]{\mathbin{\lower 4pt%
    \hbox{${\stackrel{\textstyle .}{\scriptscriptstyle #1}}$}}}

\newcommand{\du}{\varhexstar} 

\newcommand{\duer}{%
\mathbin{\raisebox{3pt}{\varhexstar}\kern-3.70pt{\rule{0.15pt}{4pt}}}\,}

\usepackage{fancyhdr}
\begin{document}

\title{{\bf Ehresmann doubles
and Drinfel'd doubles\\ 
for Lie algebroids
and Lie bialgebroids}
} 

\author{K. C. H. Mackenzie\\
        Department of Pure Mathematics\\
        University of Sheffield\\
        Sheffield, S3 7RH\\
        United Kingdom\\
        {\sf K.Mackenzie@sheffield.ac.uk}}

\date{December 18, 2006}

\maketitle

\begin{abstract}
The word `double' was used by Ehresmann to mean `an object $X$ in
the category of all $X$'. Double categories, double groupoids and 
double vector bundles are instances, but the notion of Lie 
algebroid cannot readily be doubled in the Ehresmann sense, since 
a Lie algebroid bracket cannot be defined diagrammatically. In 
this paper we use the duality of double vector bundles to define
a notion of double Lie algebroid, and we show that this abstracts 
the infinitesimal structure (at second order) of a double 
Lie groupoid. 

We further show that the cotangent of either Lie algebroid in a 
Lie bialgebroid has a double Lie algebroid structure, and that a
pair of Lie algebroid structures on dual vector bundles forms a
Lie bialgebroid if and only if the structures which they canonically
induce on their cotangents form a double Lie algebroid. In 
particular, the Drinfel'd double of a Lie bialgebra has a double 
Lie algebroid structure. 

We also show that matched pairs of Lie algebroids, as used
by J.--H. Lu in the classification of Poisson group actions, are in
bijective correspondence with vacant double Lie algebroids. 
\end{abstract}

\newpage

\section*{INTRODUCTION}

The word `double' has two distinct usages. 
In \cite{Drinfeld:1987}
Drinfel'd introduced the classical double of a Lie bialgebra
$(\gog, \gog^*)$, which combines $\gog$ and $\gog^*$ into a single
Lie algebra $\gog\bowtie\gog^*$ 
which encodes fully the relations between $\gog$ and $\gog^*$. 
(We do not consider quantum doubles in this paper.)

In categorical language the word `double' describes an object
of type $X$ in the category of all objects $X$. Thus a `double
vector bundle' is a vector bundle $D\to A$ in which both the 
total space and base space are objects in the category of vector 
bundles, say $D\to B$ and $A\to M$; the maps defining the
structure of $D\to A$ must be morphisms in the category of vector 
bundles. (More detail is given below and in \S\ref{sect:ctvb}.) 
The notion of a double category, 
the fundamental prototype of this concept of double, 
was introduced and developed by Ehresmann \cite{Ehresmann:CetS}. 
Double categories should be distinguished from the more recent 
notions of 2--category and bicategory. 

We show in this paper that the Drinfel'd double of a Lie bialgebra
may be regarded as a double in the Ehresmann sense. Further, we show
that the cotangent of (either bundle in) a Lie bialgebroid is a
double in the Ehresmann sense, and that this construction reduces 
to the Drinfel'd double when the Lie bialgebroid is a Lie bialgebra. 

The crucial step is to define a concept of double in the Ehresmann
sense for Lie algebroids. The definition cannot be precisely of the
type described above since the bracket of a Lie algebroid is not 
a map in the category of vector bundles. We use an indirect approach 
(Definition \ref{df:doubla}) and show, in \S\ref{sect:dladlg}, that
this is an infinitesimal form of doubles of Ehresmann type. 

A double Lie algebroid is first of all a double vector bundle as in
Figure \ref{fig:intro}(a);
\begin{figure}[h]
\begin{picture}(340,100)(-30,0)
\put(0,90){\xymatrix@=15mm{
D \ar[r] \ar[d] & B \ar[d] \\
A \ar[r] & M 
}}
\put(30,0){(a)}
\put(150,90){\xymatrix@=15mm{
TA \ar[r] \ar[d] & TM \ar[d] \\
A \ar[r] & M 
}}
\put(190,0){(b)}
\put(300,90){\xymatrix@=15mm{
T^*A \ar[r] \ar[d] & A^* \ar[d] \\
A \ar[r] & M 
}}
\put(340,0){(c)}
\end{picture}
\caption{\ \label{fig:intro}}
\end{figure}
that is, $D$ has two vector bundle structures, on bases $A$
and $B$, each of which is itself a vector bundle on base $M$, such 
that for each structure on $D$, the structure maps (projection, 
addition, scalar multiplication) are vector bundle morphisms with 
respect to the other structure, and the combination of the two
projections, $D\to A\times_M B$, is a surjective submersion
(see \S\ref{sect:ctvb} for more 
detail). Two examples to keep in mind are the tangent
prolongation of an ordinary vector bundle as in 
Figure~\ref{fig:intro}(b) (see \cite{Besse} for a classical 
treatment), and the cotangent double vector bundle as in 
Figure~\ref{fig:intro}(c) (see \cite{MackenzieX:1994}).

Now suppose that all four sides of Figure~\ref{fig:intro}(a) have Lie
algebroid structures. The problem is to define compatibility between 
the bracket structures on $D$, bearing in mind that the brackets
are on different modules of sections. The key is the duality for 
double vector bundles described in \S\ref{sect:ctvb}. Although our 
compatibility condition cannot (as far as we know) be interpreted 
as a condition of the form `either bracket structure on
(sections of) $D$ is a morphism with respect to the other', we
show in Theorem~\ref{thm:dladlg} that this concept of double Lie 
algebroid includes the infinitesimal invariants of double Lie 
groupoids, which \emph{are} defined by such conditions. 
We therefore feel justified in referring to it as a double of
Ehresmannian type, or a categorical double. 

In \S\ref{sect:dlalba} we define a notion of double for Lie 
bialgebroids. A Lie bialgebroid $(A,A^*)$ on base $M$ 
\cite{MackenzieX:1994}
consists of Lie algebroid structures on a vector bundle and its 
dual, subject to a compatibility condition which can be most 
succinctly expressed by saying that the differential of each
Lie algebroid acts as a derivation on the Schouten (or Gerstenhaber)
algebra associated with the other \cite{Kosmann-Schwarzbach:1995}. 

Given a Lie bialgebroid $(A,A^*)$, the Lie algebroid structure on 
$A$ induces a Poisson structure on $A^*$ by an extension 
\cite{Courant:1990t} of the duality between Lie algebras and 
linear Poisson structures on vector spaces. This Poisson structure 
in turn induces a Lie algebroid structure on
$T^*A^*\to A^*$, defined in terms of the Poisson bracket of 
1--forms. Likewise, the Lie algebroid structure on $A^*$ induces a 
Lie algebroid structure on $T^*A\to A.$ Using the canonical 
diffeomorphism $T^*A^*\to T^*A$ \cite{MackenzieX:1994} we now have
Lie algebroid structures on each of the four sides of 
Figure~\ref{fig:intro}(c). This is a double Lie algebroid as
defined in \S\ref{sect:adla}, the \emph{cotangent double of}
$(A, A^*).$ 

Suppose now that we have Lie algebroid structures on $A$ and $A^*$
with no further assumption. As above, the Poisson structures on
$A^*$ and $A$ induce Lie algebroid structures on $T^*A^*\to A^*$
and $T^*A\to A$ and we prove in \S\ref{sect:dlalba} that if these
constitute a double Lie algebroid structure on $T^*A$, then 
$(A,A^*)$ is a Lie bialgebroid. 

When $M$ is a point, so that $A = \gog$ and $A^*$ are Lie algebras, 
the cotangent manifold $T^*\gog$ is $\gog\times\gog^*$. The Lie 
algebroid structure on $\gog\times\gog^*\to\gog$
is the action, or transformation, Lie algebroid defined by the
coadjoint action of $\gog^*$ on $\gog$; likewise the coadjoint
action of $\gog$ on $\gog^*$ defines the Lie algebroid structure
on $\gog\times\gog^*\to\gog^*.$
These structures have not been much used, since $\gog\times\gog^*$
has the classical Drinfel'd double structure $\gog\bowtie\gog^*$, 
as a Lie algebra.

For an arbitrary Lie bialgebroid $(A,A^*)$ there is no known Lie
algebroid structure on $T^*A\to M$ (and indeed there is no natural
vector bundle structure) which reduces to 
$\gog\bowtie\gog^*$ in the bialgebra case.  Consequently there is
no direct analogue for $T^*A$ of the Manin triple theorem 
(as in \cite{LuW:1990})  which, for a pair of Lie algebra 
structures on dual vector 
spaces, gives a criterion for them to form a Lie bialgebra. 
Nonetheless Theorem \ref{theorem:Manin} provides a 
criterion for a pair of Lie algebroid structures on dual vector 
bundles to form a Lie bialgebroid. 

A Manin triple theorem for Lie bialgebroids which directly
generalizes the Lie bialgebra result, was given by Liu, 
Weinstein, and Xu \cite{LiuWX:1997}, whose notion of double is a 
Courant algebroid structure on the direct sum $A\oplus A^*$. Some 
comments on the differences between the two approaches are given 
at the end of \S\ref{sect:dlalba}. 

The lack of a structure on $T^*A\to M$ is an inherent feature of
the cotangent. A cotangent double vector bundle $T^*A$ is 
locally isomorphic to $A\times_M A^*\times_M T^*M$ and is of the
form $A\times_M A^*$ if and only if $M$ is $0$--dimensional. 

We call a double Lie algebroid $(D;A,B;M)$
\emph{vacant} if $D$ is isomorphic to $A\times_M B$ as a double
vector bundle. We prove in \S\ref{sect:mpvdla} that in this case 
there are induced actions $\rho$ of $A$ on $B$ and $\sigma$ of $B$ 
on $A$, and these satisfy the equations
\begin{gather}                                   
\rho_X([Y_1, Y_2]) = [\rho_X(Y_1), Y_2] + [Y_1, \rho_X(Y_2)]
    +\rho_{\sigma_{Y_2}(X)}(Y_1) - \rho_{\sigma_{Y_1}(X)}(Y_2),
\nonumber\\
\sigma_Y([X_1, X_2]) = [\sigma_Y(X_1), X_2] + [X_1, \sigma_Y(X_2)]
    +\sigma_{\rho_{X_2}(Y)}(X_1) - \sigma_{\rho_{X_1}(Y)}(X_2),
\label{eq:first}\\
a(\sigma_Y(X)) - b(\rho_X(Y)) = [b(Y), a(X)],\nonumber
\end{gather}
for $X_1, X_2, X\in\Ga A,\ Y_1, Y_2, Y\in\Ga B$, where $a$ and $b$
are the anchors. Conversely, given such a pair of actions, 
$A\times_M B$ has a double Lie algebroid structure. 

A pair of actions satisfying (\ref{eq:first}) is called a 
\emph{matched pair} structure for $A$ and $B$ \cite{Mokri:1997}.
This notion of matched pair extends the concept introduced into
geometry by Kosmann--Schwarzbach and Magri 
\cite{Kosmann-SchwarzbachM:1988} under the name 
\emph{extension bicrois\'ee} or \emph{twilled extension}, by Lu 
and Weinstein \cite{LuW:1990} as a \emph{double Lie algebra}, and 
by Majid \cite{Majid:1990}, with the term
\emph{matched pair} which we use. (In fact, forms of the concept 
had been found much earlier; see \cite{Weinstein:1990}, 
\cite{ChariP:1994} for references.) As in the case of Lie 
algebras, given a matched pair structure for Lie algebroids $A$
and $B$, the direct sum bundle $A\oplus B$ has a Lie algebroid 
structure $A\bowtie B$ for which $A\oplus 0$ and $0\oplus B$ are Lie 
subalgebroids, and the converse is also true \cite{Mokri:1997}.

We note here that Lu \cite{Lu:1997} showed that any Poisson action 
of a Poisson Lie group $G$ on a Poisson manifold $P$ gives 
rise to a matched pair of Lie algebroids, and the corresponding
diagonal structure $T^*P\bowtie(P\times\gog^*)$ carries
Drinfel'd's classification of Poisson homogeneous spaces 
\cite{Drinfeld:1993}. In a future paper we will show that the
corresponding construction for an action of a general Poisson
groupoid gives rise to a double Lie algebroid which is not
vacant, and so does not correspond to a matched pair. 

\minilos 

There are two further points to be made about the concept of
double Lie algebroid. Firstly, a natural global analogue 
is provided by the concept of double Lie groupoid 
\cite{Mackenzie:1992}: we show in \S\ref{sect:dladlg} that 
applying a second--order version of the Lie functor to a double
Lie groupoid yields a double Lie algebroid. This subsumes various
known differentiation results and raises a general integrability
problem. 

It is well known that Lie bialgebroids are the infinitesimal form
of Poisson groupoids \cite{Weinstein:1988}, \cite{MackenzieX:2000}. 
In the case of Lie bialgebras a further step was achieved by Lu
and Weinstein \cite{LuW:1989} who proved the existence of a
symplectic double Lie groupoid integrating a Poisson Lie group. 
The double Lie algebroid of this symplectic double Lie groupoid
is the cotangent double Lie algebroid described above. 

For matched pairs of Lie algebras the global concept of
matched pair of Lie groups was introduced into geometry at the same
time \cite{Majid:1990}, \cite{LuW:1990}. The global concept of
a matched pair of Lie groupoids was given in \cite{Mackenzie:1992}. 
In the Lie algebra case, various integrability results were given in
the same papers, but little is known about the integrability
of matched pairs of Lie algebroids \cite{MackenzieM:1999}. As
\cite{LuW:1989} shows, a matched pair of Lie algebras may fail to
integrate to a matched pair of Lie groups, but nonetheless
integrate to a double Lie groupoid. 

The second point to be made is that the concept of double Lie
algebroid offers, at least in some situations, a means of bypassing
the problems caused by the lack of pullbacks in the Poisson
category. It is well--known that the action groupoid arising from
the action of a Poisson Lie group $G$ on a Poisson manifold $P$ is 
not a Poisson groupoid. It was shown in \cite{Mackenzie:2000lmp} that
the induced action of $T^*G\gpd \gog^*$ on the dual $\gop\co
T^*P\to\gog^*$ of the infinitesimal action, called the \emph{pith} 
in \cite{Mackenzie:2000lmp}, gives rise to a groupoid object in
the category of Lie algebroids which differentiates to Lu's matched 
pair of Lie algebroids. In the group case one may work directly with 
the matched pair but, as mentioned above, for a Poisson action of a 
Poisson groupoid, the corresponding structure is a double Lie 
algebroid which does not reduce to a matched pair. 

Since the mid 1990s several authors have worked to provide a
characterization of the cotangent double of a Lie bialgebroid in
terms of super geometry \cite{RoytenbergW:1998}, 
\cite{Roytenberg:thesis}. 
This has now been achieved by Ted Voronov \cite{Voronov:mtqm} who
further gives a formulation of the general concept of double Lie
algebroid presented here in terms of super geometry. 
Independently, Mehta \cite{Mehta:thesis} has found a super 
formulation of the related notion of \LAgpd. 

On a speculative level, one may ask what relation, if any, exists
between this theory and the weak multiple category theory that
is used in some approaches to 
field theory. The differences between double Lie groupoid theory 
\cite{Mackenzie:1992}, \cite{Mackenzie:2000} 
and the various notions of weak 2--category or bicategory, are 
very considerable and it may be that there is no simple 
relationship. The expectation that a link may exist is perhaps 
based on 
the fact that the notion of ordinary Lie groupoid includes both 
the frame groupoids of differential geometry and the fundamental 
groupoid of a manifold, and thus has both a differential geometric 
and a topological significance. However, this may be misleading 
as regards multiple structures: the homotopy properties of 
double Lie groupoids present considerable problems, and the notion 
of weak 2--category does not seem to readily accept a smooth 
structure. Perhaps these matters can be resolved, 
but in any case the results of this paper, particularly in the light
of \S\ref{sect:dladlg}, appear to show that multiple Lie groupoid 
theory is well suited to finite--dimensional smooth differential 
geometry. 

\minilos

An announcement of some of the main results of this paper was 
published in \cite{Mackenzie:1998}, and an early version was
arXived in \cite{Mackenzie:notions}. Amongst the many changes 
between \cite{Mackenzie:notions} and the current paper, the most
striking is perhaps the simplification of the definition 
\ref{df:doubla} of a double Lie algebroid, following the observation
of Ted Voronov that in his super formulation \cite{Voronov:mtqm}, 
the analogue of Condition II in \cite{Mackenzie:notions} follows
from Conditions I and III. This has enabled many arguments to be
considerably shortened. 

We do not include coordinate representations in this paper. The use 
of intrinsic formulations keeps clear the nature of each structure 
and makes the lifting of results and definitions from one level to
another much easier. For example, most of the canonical isomorphisms 
used in the paper have been introduced for manifolds or vector 
bundles, and are then lifted to various multiple structures. The
functorial nature of these processes is less clear in 
coordinate formulations. However, a coordinate--based approach is 
included in \cite{Voronov:mtqm}. 

The material of this paper has developed over a considerable
number of years and I have benefited greatly from many conversations
with Yvette Kos\-mann--Schwarz\-bach, Alan Weinstein, and Ping Xu
over this time. Conversations with Ted Voronov at the 
\emph{XXVth Workshop on Geometry and Mathematical Physics} at 
Bia{\l}owie{\.z}a in 2006 provided a crucial stimulus for 
this revision and I am very grateful to him, and to the organizers 
of the Bia{\l}owie{\.z}a workshop for providing, 
inter alia, the opportunity for extended conversations. 
I also wish to thank the organizers of the programme on 
\emph{Gerbes, Groupoids and Quantum Field Theory} at the ESI in 
2006, where important parts of this revision were written. 

\newpage

\section{DOUBLE VECTOR BUNDLES AND THEIR DUALITY}
\label{sect:ctvb}

We begin by recalling the duality of double vector bundles which was 
introduced in \cite{Mackenzie:1999}, \cite{KoniecznaU:1999} using 
the duality which Pradines \cite{Pradines:1988} had introduced for 
vector bundle objects in the category of Lie groupoids. See 
\cite{Mackenzie:GT} for a fuller account. 

The definition of \emph{double vector bundle} given in the 
Introduction (see
Figure~\ref{fig:intro}(a)) is complete but needs elaboration. 
We denote the projections in Figure~\ref{fig:intro}(a) by
$q^D_A\co D\to A$ and $q^D_B \co D\to B$; these are vector 
bundle morphisms over the projections $q_B\co B\to M$ and 
$q_A\co A\to M$ respectively. If $d\in D$ has $q^D_A(d) = a$ and 
$q^D_B(d) = b$, and $q_A(a) = m$, then we say that $d$ has 
\emph{outline} $(d;a,b;m).$

Consider four elements $d_1,\dots,d_4$ of $D$ with projections 
$a_1,\dots,a_4$ in $A$ and $b_1,\dots,b_4$ in $B$ such that
$a_1 = a_2$, $a_3 = a_4$, $b_1 = b_3$ and $b_2 = b_4$. Then the 
condition that each bundle projection from $D$ is a morphism of 
vector bundles with respect
to the other structure guarantees that both sides of 
\begin{equation}
\label{eq:ichange}
(d_1\add{A} d_2)\add{B}(d_3\add{A} d_4) =
(d_1\add{B} d_3)\add{A}(d_2\add{B} d_4)
\end{equation}
are defined. Here $\add{A}$ denotes the addition in $D\to A$ 
and $\add{B}$ the addition in $D\to B$. 
The condition that each addition be a vector bundle morphism with 
respect to the other structure then implies equality in 
(\ref{eq:ichange}). There 
are corresponding equations involving scalar multiplication. Equation
(\ref{eq:ichange}) is the \emph{interchange law for addition}. 

Given $a\in A$ we denote the zero of $D\to A$ above $a$
by $\tilzero_a$. Similarly the zero of $D\to B$ above
$b\in B$ is denoted $\tilzero_b$. The zeros in $A\to M$ and $B\to M$
are denoted $0^A_m$ and $0^B_m$. The interchange law implies that 
the zero $\tilzero_{0^A_m}$ for $D\to A$ and the zero 
$\tilzero_{0^B_m}$ for $D\to B$ are equal; this is denoted 
$\odot_m$. 

The intersection of the kernels of
the two projections defined on $D$ is called the \emph{core} of
the double vector bundle $D$ and denoted $C$. The vector bundle
structures on $D$ induce a common vector bundle structure on $C$
with base $M$. The kernel of $q^D_B\co D\to B$ is the
pullback $q_A\xclam C$ of this vector bundle $C$ across $q_A$, and 
the kernel 
of $q^D_A\co D\to A$ is $q_B\xclam C$. In practice cores are
often identified with structures which are not strictly subsets of
$D$ and in such cases we write $\Bar{c}$ for the element of
$D$ which corresponds to a $c\in C$.

A \emph{morphism of double vector bundles from $(D;A,B;M)$
to $(D';A',B';M')$} is a quadruple of maps $\phi\co D\to D'$, 
$\phi_A\co A\to A'$, $\phi_B\co B\to B'$ and $\phi_M\co M\to M'$
such that $(\phi, \phi_A)$, $(\phi, \phi_B)$, $(\phi_A, \phi_M)$
and $(\phi_B, \phi_M)$ are morphisms of vector bundles. A morphism
of double vector bundles induces a morphism of the cores
$\phi_C\co C\to C'$ by restriction. 

Regarding $D$ as an ordinary vector bundle over $A$, we denote the 
dual
bundle by $D\duer A$. (This notation avoids problems with multiple 
subscripts later.) In addition to its standard structure on base $A$,
this has a vector bundle structure on base $C^*$, the projection of 
which, denoted $\gamma^A_{C^*}$, is defined by 
\begin{equation}                             
\label{eq:unfproj}
\langle\gamma^A_{C^*}(\Phi),\, c\rangle =
      \langle\Phi,\, \tilzero_a \add{B} \Bar{c}\rangle
\end{equation}
where $\Phi\co (q^D_A)^{-1}(a)\to\R,\ a\in A_m,$ and $c\in C_m$.
The addition in $D\duer A\to C^*$, which we denote by $\add{C^*}$,
is defined by
\begin{equation}                  
\langle\Phi\add{C^*}\Phi',\, d \add{B} d'\rangle =
   \langle\Phi,\,d\rangle + \langle\Phi',\,d'\rangle.
\end{equation}
The condition that $\Phi$ and $\Phi'$ project to the same element 
of $C^*$ guarantees that this is well--defined. The zero of 
$D\duer A\to C^*$ above $\kappa\in C^*_m$ is denoted 
$\tilzero^{\du A}_\kappa$ and is defined by
\begin{equation}
\label{eq:dualzero}
\langle\,\tilzero^{\du A}_\kappa,\, \tilzero_b \add{A} \Bar c\rangle 
= \langle\kappa,\, c\rangle
\end{equation}
where $c\in C_m$ and $b$ is any element of $B_m$. The scalar 
multiplication is defined in a similar way. These two structures 
make $D\duer A$ into a double
vector bundle as in Figure~\ref{fig:dvbduals}(a),
\begin{figure}[h]
\begin{picture}(340,100)(-30,0)
\put(50,90){\xymatrix@=15mm{
D\duer A \ar[r]^{\gamma^A_{C^*}} \ar[d]_{\gamma^A_A} & 
         C^* \ar[d]^{q_{C^*}} \\
A \ar[r]_{q_A} & M 
}}
\put(95,0){(a)}
\put(220,90){\xymatrix@=15mm{
D\duer B \ar[r]^{\gamma^B_B} \ar[d]_{\gamma^B_{C^*}} & 
          B \ar[d]^{q_B} \\
C^* \ar[r]_{q_{C^*}} & M 
}}
\put(265,0){(b)}
\end{picture}
\caption{\ \label{fig:dvbduals}}
\end{figure}
the \emph{vertical dual of $D$}. The core of $D\duer A$
identifies with $B^*$, with the core element $\Bar\psi$
corresponding to $\psi\in B_m^*$ given by
\begin{equation}
\label{eq:corepsi}
\langle\Bar\psi,\, \tilzero_b \add{A} \Bar c\rangle =
\langle\psi,\,b\rangle.
\end{equation}

There is also a \emph{horizontal dual} $D\duer B$ with sides $B$
and $C^*$ and core $A^*$, as in Figure~\ref{fig:dvbduals}(b). 
There is
now the following somewhat surprising result. A proof is in
\cite[Chap.~9]{Mackenzie:GT}. 

\begin{thm}[{\cite[3.1]{Mackenzie:1999}}, 
{\cite[Thm.~16]{KoniecznaU:1999}}]
\label{thm:dualduality} 
There is a natural (up to sign) nondegenerate pairing between 
$D\duer A\to C^*$ and $D\duer B\to C^*$ given by
\begin{equation}                       \label{eq:3duals}
\thbr{\Phi}{\Psi} = \langle\Phi,\, d \rangle 
                                  - \langle d,\,\Psi \rangle
\end{equation}
where $\Phi\in D\duer A,\ \Psi\in D\duer B$ have
$\gamma^A_{C^*}(\Phi) = \gamma^B_{C^*}(\Psi)$ and $d$ is any element
of $D$ with $q^D_A(d) = \gamma^A_A(\Phi)$ and
$q^D_B(d) = \gamma^B_B(\Psi).$
\end{thm}

The pairing on the LHS of (\ref{eq:3duals}) is over $C^*$, whereas 
the pairings on the RHS are over $A$ and $B$ respectively. The 
pairing induces isomorphisms of double vector bundles
\begin{eqnarray}
Z_A\co D\duer A \to D\duer B\duer C^*, & \quad & 
\langle Z_A(\Phi),\,\Psi\rangle = \thbr{\Phi}{\Psi},\\
Z_B\co D\duer B \to D\duer A\duer C^*, & \quad & 
\langle Z_B(\Psi),\,\Phi\rangle = \thbr{\Phi}{\Psi}
\end{eqnarray}
with $(Z_A)\duer C^* = Z_B.$ Note that $Z_B$ induces the identity on 
both
side bundles $C^*$ and $B$ and is $-\id$ on the cores $A^*$, whereas 
$Z_A$ induces the identity on the side bundles $C^*$ and the cores 
$B^*$ but is $-\id$ on the side bundles $A.$

The example which follows is needed for \S\ref{sect:dlalba}. 

\begin{ex}\rm
\label{ex:tangentd}
Given an ordinary vector bundle $(A,q,M)$ there is the tangent double
vector bundle of Figure~\ref{fig:intro}(b). This example was used 
as early as the 1960s to provide an efficient formalism for 
connection theory in $A$; see \cite{Besse}. The core of $TA$ consists
of the vertical vectors along the zero section and identifies 
canonically with $A$. Dualizing the standard structure $TA\to A$ 
gives the cotangent double vector bundle of 
Figure~\ref{fig:intro}(c) with core $T^*M$. Dualizing
the structure $TA\to TM$ gives a double vector bundle which we denote
$(T^\sol A;A^*,TM;M)$; this is canonically isomorphic to 
$(T(A^*);A^*,TM;M)$ under an isomorphism $I_A\co T(A^*)\to T^\sol A$ 
given by
\begin{equation}
\label{eq:I}
\langle I_A(\mathcal{X}),\,\xi \rangle =
\llangle\mathcal{X},\, \xi \rrangle
\end{equation}
where $\mathcal{X}\in T(A^*)$ and $\xi \in TA$ have
$T(q_*)(\mathcal{X}) = T(q)(\xi)$, and $\llangle\ ,\ \rrangle$ is the
tangent pairing of $T(A^*)$ and $TA$ over $TM$ defined by
$$
\llangle\mathcal{X},\ \xi\rrangle = 
\left.\tfrac{d}{dt}\langle \phi_t,a_t\rangle\right|_0,
$$
for $\mathcal{X} = \left.\frac{d}{dt} \phi_t\right|_0$ and 
$\xi = \left.\frac{d}{dt} a_t\right|_0$. 

Applying \ref{thm:dualduality} to $D = TA$ thus induces a 
pairing of $T^*A$ and $T(A^*)$ over $A^*$ given by
$$
\thbr{\Phi}{\mathcal{X}} =
\llangle\mathcal{X},\xi\rrangle -
\langle\Phi, \xi\rangle
$$
where $\Phi\in T^*A$ and $\mathcal{X}\in T(A^*)$ have
$r(\Phi) = p_{A^*}(\mathcal{X})$, and $\xi\in TA$ is chosen so that
$T(q)(\xi) = T(q_*)(\mathcal{X})$ and $p_A(\xi) = c_A(\Phi)$.
(Here $c_A$ is the projection of the standard cotangent bundle and
$r\colon T^*A\to A^*$ is a particular case of (\ref{eq:unfproj}).)
This pairing is nondegenerate and so defines an isomorphism of 
double vector bundles $R\co T^*A^*\to T^*A$ by the condition
$$
\thbr{R(\mathcal{F})}{\mathcal{X}} =
            \langle\mathcal{F},\, \mathcal{X}\rangle
$$
where the pairing on the RHS is the standard one of $T^*(A^*)$ and
$T(A^*)$ over $A^*$. This $R$ preserves the side bundles $A$ and 
$A^*$ but induces $-\id\co T^*M\to T^*M$ as the map of cores. 
It is an antisymplectomorphism with respect to the exact 
symplectic structures. In summary we now have the very useful 
equation
\begin{equation}                                   
\label{eq:vue}
\langle\mathcal{F},\, \mathcal{X}\rangle + \langle R(\mathcal{F}),\, \xi\rangle
= \llangle\mathcal{X},\, \xi\rrangle,
\end{equation}
for $\mathcal{F}\in T^*A^*,\ \mathcal{X}\in T(A^*),\ \xi\in TA$, 
where the pairings are over $A^*, A$ and $TM$ respectively. This 
equation was originally derived in \cite[\S5]{MackenzieX:1994} 
by using the symplectic structures. 
\end{ex}

Now return to the general double vector bundle in 
Figure~\ref{fig:intro}(a), denoting the core by $C$. Each vector 
bundle structure on $D$ gives rise to a double cotangent bundle. 
These two double cotangents fit together into a triple structure 
as in Figure~\ref{fig:cottrip}(a). This structure is a triple 
vector bundle in an obvious sense: each face is a double vector 
bundle and for each vector bundle structure on $T^*D$, the maps 
defining the structure are morphisms of double vector bundles. The 
left and rear faces of Figure~\ref{fig:cottrip}(a) are the two 
cotangent doubles of the two structures on $D$ and the top face 
may be regarded as the cotangent double of either of the two
duals of $D$. (In all such triple diagrams, we take the oblique
arrows to be coming out of the page.) 

\begin{figure}[h]
\begin{picture}(340,150)(0,30)
\put(0,150){\xymatrix@=5mm{
T^*D \ar[rr] \ar[dd] \ar[dr] 
& & D\duer B \ar[dd] \ar[dr] &\\
& D\duer A  \ar[rr] \ar[dd]         & & C^* \ar[dd]  \\
D    \ar[rr] \ar[dr]         & & B \ar[dr] &\\
& A  \ar[rr] & & M\\
}}
\put(100,30){(a)}
\put(250,150){\xymatrix@=6mm{
TD \ar[rr] \ar[dd] \ar[dr] 
& & TB \ar[dd] \ar[dr] &\\
& TA  \ar[rr] \ar[dd]         & & TM \ar[dd]  \\
D    \ar[rr] \ar[dr]         & & B \ar[dr] &\\
& A  \ar[rr] & & M\\
}}
\put(320,30){(b)}
\end{picture}
\caption{\ \label{fig:cottrip}}
\end{figure}

Figure~\ref{fig:cottrip}(a) is, in a sense made precise in 
\cite{Mackenzie:2005dts}, the vertical dual of the tangent 
prolongation of $D$ as given in Figure~\ref{fig:cottrip}(b). 
Five of the six faces in \ref{fig:cottrip}(a) are double vector 
bundles of types considered already; it is only necessary
to verify that the top face is a double vector bundle. The cores 
of the same five faces are known, and we take the core of the top 
face to be $T^*C$.
Taking these cores in opposite pairs, together with a parallel 
edge, then gives three double vector bundles: the left and the 
right cores form $(T^*A;A,A^*;M)$, the back and the front cores 
form $(T^*B;B,B^*;M)$ and the up and the down cores form 
$(T^*C;C,C^*;M)$. Each of these \emph{core double vector bundles} 
has core $T^*M$.

\newpage

\section{PRELIMINARY CASE}
\label{sect:pc}

Consider a double vector bundle as in 
Figure~{\rm \ref{fig:intro}(a)}. Before considering bracket 
structures, we need to consider two classes of sections of the 
bundle structures on $D$. 

A section $\xi\in\Ga_A D$ is a \emph{(vertical) linear section}
if there exists $X\in\Ga B$ such that $(\xi, X)$ is a vector 
bundle morphism from $(A, q_A, M)$ to $(D, q^D_B, B)$. 

\begin{prop}
If $(\xi, X)$ is a linear section, then $\ell_\xi\co D\duer A\to\R$ 
defined by
$$
\Phi\mapsto \langle \Phi,\,\xi(\gamma^A_A(\Phi))\rangle
$$
is linear with respect to $C^*$ as well as $A$, and the restriction 
of $\ell_\xi$ to the 
core of $D\duer A$ is $\ell_X\co B^*\to\R$.
\end{prop}

\pf
Take elements $(\Phi_1; a_1, \kappa; m)$ and 
$(\Phi_2; a_2, \kappa; m)$ of $D\duer A$. Their sum over $C^*$ 
has outline
$(\Phi_1\add{C^*}\Phi_2; a_1 + a_2, \kappa; m)$ and so, using 
first the
linearity of $\xi$ and then the definition of $\add{C^*}$,
\begin{eqnarray*}
\ell_\xi(\Phi_1\add{C^*}\Phi_2)
& = & \langle \Phi_1\add{C^*}\Phi_2,\, \xi(a_1 + a_2)\rangle =
\langle \Phi_1\add{C^*}\Phi_2,\, \xi(a_1) \add{B} \xi(a_2)\rangle\\
& = & 
\langle \Phi_1, \xi(a_1)\rangle + \langle \Phi_2, \xi(a_2)\rangle 
= \ell_\xi(\Phi_1) + \ell_\xi(\Phi_2).
\end{eqnarray*}
Scalar multiplication is similar. Given $\psi\in B^*_m$ the 
corresponding core element $\Bar{\psi}\in D\duer A$ is given by 
(\ref{eq:corepsi}). Using this, 
\begin{equation}
\label{eq:xix2}
\ell_\xi(\Bar{\psi}) = 
\langle \Bar{\psi},\, \xi(0^A_m) \rangle 
= \langle \psi,\, x(m) \rangle = \ell_x(\psi),
\end{equation}
where $\ell_x\co B^*\to\R$ denotes the linear
map $\phi\mapsto \langle \phi,\,x(q_{B^*}(\phi))\rangle,$
for any $x\in\Ga B$. 
\pfend

Given a linear section $(\xi, x)$, denote by $\xi^\sqcap$ the section
of $D\duer A\duer C^*\to C^*$ which is defined by $\ell_\xi$. Thus
\begin{equation}
\label{eq:xisqcap}
\langle{\xi^\sqcap(\kappa)},\,\Phi\rangle_{C^*}
= \ell_\xi(\Phi)
= \langle\Phi,\, \xi(a)\rangle_A
\end{equation}
for $(\Phi;a, \kappa; m)\in D\duer A$. It is easily seen that
$(\xi^\sqcap, x)$ is a linear section of $D\duer A\duer C^*\to C^*$.

There is thus a bijective correspondence between vertical linear 
sections of $D$ and horizontal linear sections of 
$D\duer A\duer C^*$. In the case where $D$ is the tangent of an 
ordinary vector bundle $A$, as in Figure~\ref{fig:intro}(b),
this is the correspondence between linear vector fields on $A$ 
and linear vector fields on $A^*$ (see \cite[\S3.4]{Mackenzie:GT}).

Next, for any $c\in\Ga C$, define a \emph{core section}
$\Bar{c}\in\Ga_A D$ by
$$
\Bar{c}(a) = \Bar{c(m)}\add{B}\tilzero_a.
$$

\begin{prop}
\label{prop:gen}
The vertical linear sections $(\xi,X)$ and the core sections
$\Bar{c}$ for $c\in\Ga C$, generate $\Ga_A D$ as a
$\cinf{A}$--module. 
\end{prop}

\pf
Firstly, given any $X\in\Ga B$, there is a linear section $(\xi,X)$. 
To see this, let 
$$
S\co D\to A\times_M B\times_M C
$$ 
be a splitting of $D$; that is, an isomorphism of double vector 
bundles to the double vector bundle for which the vertical structure
is the pullback $q_A^!(B\oplus C)$ and the horizontal structure
is the pullback $q_B^!(A\oplus C)$. Then the section
$$
a\mapsto S^{-1}(a,\,X(m),\,0_m), 
$$
where $m = q_A(a)$, is linear over $X$. 

Now consider any section whatever $\xi\in\Ga_A D$. Then by
paracompactness $q^D_A\circ\xi\co A\to B$ is a finite sum 
$$
q^D_A\circ\xi = \sum F_i \Bar{X_i}
$$
where $F_i\in\cinf{A}$ and $\Bar{X_i}$, for $X_i\in\Ga B$, 
denotes the pullback section of $q_A^!B.$ For each $X_i$ choose
a linear section $(\xi_i,X_i)$ and write $\eta = \sum F_i\xi_i.$

Then $\xi - \eta$ projects to zero under $q^D_A$ and is therefore
$\Bar{c}$ for some $c\in\Ga C.$
\pfend

Before addressing the general concept of double Lie algebroid, 
it will be helpful to deal with a very special case.

\begin{df}
An \emph{\LAvb} is a double vector bundle as in 
Figure~{\rm \ref{fig:intro}(a)} together with Lie algebroid 
structures on a pair of parallel sides, such that the structure 
maps of the other pair of vector bundle structures are
Lie algebroid morphisms.
\end{df}

For definiteness, take the Lie algebroid structures to be on
$D\to A$ and $B\to M$. Write $\dan_A\co D\to TA$ and 
$b\co B\to TM$ for the anchors.

Since $D\to A$ is a Lie algebroid, the dual $D\duer A$ has
a dual Poisson structure, which is linear with respect to the bundle
structure over $A$. 

\begin{thm}                                     
\label{theorem:bothduals}
The Poisson structure $\pi_{\du A}$ on $D\duer A$ is also 
linear with respect to the bundle structure over $C^*$.
\end{thm}

\pf
Since it is dual to a Lie algebroid structure, $\pi_{\du A}$ is 
characterized by the following equations.
\begin{equation}
\label{eq:pi}
\{\ell_\xi,\,\ell_\eta\} = \ell_{[\xi,\eta]},\qquad
\{\ell_\xi,\,F\circ\gamma^A_A\} = 
\dan_A(\xi)(F)\circ\gamma^A_A,\qquad
\{F_1\circ\gamma_A^A,\,F_2\circ\gamma^A_A\} = 0,
\end{equation}
where $\xi$ and $\eta$ are arbitrary elements of $\Ga_AD$ and 
$F_1,F_2\in\cinf{A}.$ See \cite[\S10.3]{Mackenzie:GT}. We need 
to prove corresponding equations with respect to the bundle 
structure over $C^*.$ The functions $D\duer A\to\R$ which are 
linear with respect to the bundle structure over $C^*$ are 
determined by arbitrary sections of $D\duer A\duer C^*\to C^*$. 
As in Proposition \ref{prop:gen}, these are of two types. 

Firstly, if $\mscr{X},\mscr{Y}\in\Ga_{C^*}(D\duer A\duer C^*)$ 
are linear sections, then $\mscr{X} = \xi^\sqcap$ and
$\mscr{Y} = \eta^\sqcap$ for linear sections 
$(\xi,x),\,(\eta,y)\in\Ga_A D.$

\begin{lem}
\label{lem:twolinear}
If $(\xi, x)$ and $(\eta, y)$ are linear sections of $D$, then
$([\xi, \eta], [x,y])$ is also.
\end{lem}

\pf
For a section $\xi$ which projects under $q^D_B$ to a section $x$,
define a section $\xi\oplus\xi$ of $D\Whitney{B}D\to
A\Whitney{M}A$ by 
$(\xi\oplus\xi)(a_1\oplus a_2) = \xi(a_1)\oplus\xi(a_2)$.
Then $(\xi, x)$ is linear if and only if $\xi\oplus\xi$ projects to
$\xi$ under the horizontal addition $\add{B}\co D\Whitney{B}D\to D$,
which is a Lie algebroid morphism over $+\co A\Whitney{M}A\to A$
by hypothesis. Since the two components of a $\xi\oplus\xi$ depend 
on each variable separately, we have
$[\xi\oplus\xi,\, \eta\oplus\eta] = [\xi, \eta]\oplus[\xi, \eta]$ 
and the result follows.
\pfend

So $\{\ell^{C^*}_\mscr{X},\,\ell^{C^*}_\mscr{Y}\} = 
\{\ell_\xi,\,\ell_\eta\} = \ell_{[\xi,\eta]} = 
\ell^{C^*}_{[\xi,\eta]^\sqcap}$ is also linear over $C^*$. 
Here we wrote $\ell^{C^*}_\mscr{X}$ for clarity. 

Secondly, any $\phi\in\Ga A^*$ induces a \emph{core section} 
$\Bar{\phi}$ of $D\duer A\duer C^* \to C^*$ by
$$
\Bar{\phi}(\kappa) 
= \tilzero^{\du B}_\kappa \add{B} \Bar{\phi(m)},\qquad
\kappa\in C^*_m,
$$
which in turn induces $\ell^{C^*}_{\Bar{\phi}}\co D\duer A\to\R$. 
A simple calculation shows that
$$
\ell^{C^*}_{\Bar{\phi}} = \ell_\phi\circ\gamma^A_A. 
$$
It now follows from (\ref{eq:pi}) that 
\begin{equation}
\label{eq:twocores}
\{\ell^{C^*}_{\Bar{\phi}_1},\,\ell^{C^*}_{\Bar{\phi}_2}\}
= \{\ell_{\phi_1}\circ\gamma^A_A,\,\ell_{\phi_2}\circ\gamma^A_A\}
= 0.
\end{equation}
Next consider a linear section $\mscr{X} = \xi^\sqcap$ and
$\phi\in\Ga A^*$. We have
\begin{equation}
\label{eq:mixed}
\{\ell^{C^*}_\mscr{X},\,\ell^{C^*}_{\Bar{\phi}}\} = 
\{\ell_\xi,\, \ell_\phi\circ\gamma^A_A\} =
\dan_A(\xi)(\ell_\phi)\circ\gamma^A_A.
\end{equation}
Since $\dan_A(\xi)$ is a linear vector field on $A$, it maps 
the linear function $\ell_\phi$ to a linear function. 

To make this explicit, recall the notion of derivation on a vector 
bundle. A \emph{derivation} on $E\to M$ is a linear differential 
operator $\Lambda$ of order $\leq 1$ for which there is a vector 
field $X$, depending on $\Lambda$, such that 
$\Lambda(f\mu) = f\Lambda(\mu) + X(f)\mu$ for all $f\in \cinf{M}$ 
and all sections $\mu$ of $E$. For any vector bundle $E$ there is a
vector bundle $\mscr{D}(E)$ the sections of which are the 
derivations, and 
with the anchor defined by $\Lambda\mapsto X$ and the usual
bracket, $\mscr{D}(E)$ is a Lie algebroid. 
See \cite[\S3.3]{Mackenzie:GT} and references given there. 

There is a bijective correspondence between derivations on $E$ and 
linear vector fields on $E^*$, and in the present case we can 
define a derivation $\Lambda_\xi$ on $A^*$ such that
\begin{equation}
\label{eq:Lambda}
\dan_A(\xi)(\ell_\phi) = \ell_{\Lambda_\xi(\phi)}
\end{equation}
for all $\phi\in\Ga A^*.$ 

From \ref{prop:gen} we know that $\ell^{C^*}_\mscr{X}$, for 
$\mscr{X}$ an arbitrary section of $D\duer A\duer C^*\to C^*$, 
is of the form
$$
\sum F_i\,\ell^{C^*}_{\xi_i^\sqcap} + \ell^{C^*}_{\Bar{\phi}}
$$
where $F_i\in\cinf{C^*}.$ It follows from (\ref{eq:twocores}), 
(\ref{eq:mixed}), and the argument following Lemma
\ref{lem:twolinear}, that the bracket of two such functions is
another of the same type. This completes the proof that the bracket 
of any two functions $D\duer A\duer C^*\to\R$ which are linear over 
$C^*$, is also linear over $C^*$. 

We next need to consider functions which are pullbacks across
$\gamma^A_{C^*}\co D\duer A\to C^*$. Functions $C^*\to\R$ are 
of two types. Firstly any $c\in\Ga C$ induces 
$\ell_c\co C^*\to\R$. Secondly, any $f\in\cinf{M}$ pulls back 
to $f\circ q_{C^*}\co C^*\to\R$. It is easily seen that 
\begin{equation}
\label{eq:pbacks}
\ell_c\circ\gamma^A_{C^*} = \ell_{\Bar{c}},\qquad
f\circ q_{C^*}\circ\gamma^A_{C^*} = f\circ q_A\circ\gamma^A_A. 
\end{equation}

For the bracket of two functions of the first type, we have
$\{\ell_{\Bar{c}_1},\,\ell_{\Bar{c}_2}\} 
= \ell_{[\Bar{c}_1,\Bar{c}_2]}$ from (\ref{eq:pi}). Consider 
the horizontal scalar multiplication in $D$ by some $t\neq 0,1$. 
This defines an automorphism $D\to D$ over $A\to A$ and therefore 
induces a map of sections $t_B\co \Ga_{A}D \to \Ga_{A}D$; the Lie
algebroid condition then ensures that
$[t_B(\xi),\, t_B(\eta)] = t_B([\xi,\eta])$
for all $\xi, \eta\in\Ga_{A}D$. Now for $c\in\Ga C$,
$\Bar{tc} = t_B(\Bar{c})$ and so we have 
$$
\Bar{t[c_1,c_2]} = [t_B(\Bar{c}_1), t_B(\Bar{c}_2)]
= [\Bar{tc_1}, \Bar{tc_2}] = \Bar{t^2[c_1,c_2]}.
$$ 
Thus the bracket of two functions $D\duer A\to\R$ arising from
sections of $C$ must be zero. 

That the bracket of two functions of the second type in 
(\ref{eq:pbacks}) is zero follows immediately from (\ref{eq:pi}). 

Again from (\ref{eq:pi}), we have 
$$
\{\ell_c\circ\gamma^A_{C^*},\,f\circ q_{C^*}\circ\gamma^A_{C^*}\} = 
\{\ell_{\Bar{c}},\,f\circ q_{C^*}\circ\gamma^A_{C^*}\} = 
\dan_A(\Bar{c})(f\circ q_{C^*})\circ\gamma^A_{C^*}. 
$$
The anchor $\dan_A\co D\to TA$ is a morphism of double vector 
bundles and therefore restricts to a map of the cores, which we
denote by $\da_A\co C\to A$. Then $\dan_A$ maps the core section 
$\Bar{c}$ to the core section of $TA$ corresponding to $\da_A(c)$. 
The core sections of $TA$ are the vertical vector fields on $A$
and vertical vector fields map pullback functions to zero. (See
\cite[\S3.4]{Mackenzie:GT}.) So
$\{\ell_c\circ\gamma^A_{C^*},\,f\circ q_{C^*}\circ\gamma^A_{C^*}\} 
= 0.$

It remains to show that brackets
$\{\ell^{C^*}_\mscr{X},\,G\circ\gamma^A_{C^*}\}$, where
$\mscr{X}$ is a section of $D\duer A\duer C^*\to C^*$ and
$G\in\cinf{C^*},$ are of the form $G'\circ\gamma^A_{C^*}$. 
As before, it suffices to consider the cases where 
$\mscr{X}$ is $\xi^\sqcap$ or $\Bar{\phi}$ and $G$ is $\ell_c$
or $f\circ q_{C^*}$. 

Firstly, we have 
$\{\ell^{C^*}_{\ell^\sqcap},\,\ell_c\circ \gamma^A_{C^*}\} = 
\{\ell_\xi,\, \ell_{\Bar{c}}\} =
\ell_{[\xi,\, \Bar{c}]}$. Now $\xi$ projects to $x$ under $q^D_B$ 
and $\Bar{c}$ projects to 0, so $[\xi, \Bar{c}]$ also projects 
to 0. It is therefore $\Bar{c'}$
for some $c'\in\Ga C$ which we denote $Q_\xi(c)$; it is easily
checked that $Q_\xi\co\Ga C\to\Ga C$ is a derivation on $C$. 

Secondly, 
$$
\{\ell^{C^*}_{\ell^\sqcap},\,f\circ q_{C^*}\circ\gamma^A_{C^*}\} =
\{\ell_\xi,\, f\circ q_A\circ\gamma^A_A\} =
\dan_A(\xi)(f\circ q_A)\circ\gamma^A_A =
x(f)\circ q_A\circ\gamma^A_A = x(f)\circ q_{C^*}\circ\gamma^A_{C^*}
$$
which is of the required form. 

In the third case we have
$$
\{\ell_\phi\circ\gamma^A_A,\, \ell_{\Bar{c}}\} =
   -\dan_A(\Bar{c})(\ell_\phi)\circ\gamma^A_A =
   -\da_A(c)^\upa(\ell_\phi)\circ\gamma^A_A,
$$
where $\da_A(c)^\upa$ is the vertical vector field on $A$ 
corresponding to the section $\da_A(c).$ Since the flows of 
$\da_A(c)^\upa$ are by scalar multiplication, we have
$\da_A(c)^\upa(\ell_\phi) = \langle\phi,\,\da_A(c)\rangle\circ q_A$ 
and the result is of the required form. 

Lastly, in the fourth case the bracket is zero, since both functions
are pullbacks from $A$. 

This completes the proof of Theorem \ref{theorem:bothduals}.
\pfend

It follows that $D\duer A\duer C^*\to C^*$ has a Lie algebroid
structure defined by 
$$ 
\ell^{C^*}_{[\mscr{X},\mscr{Y}]} = 
\{\ell^{C^*}_\mscr{X},\,\ell^{C^*}_\mscr{Y}\},\qquad
e(\mscr{X})(G)\circ\gamma^A_{C^*} = 
\{\ell^{C^*}_\mscr{X},\,G\circ\gamma^A_{C^*}\}, 
$$ 
where $e$ denotes the anchor, 
$\mscr{X},\mscr{Y}\in\Ga_{C^*}(D\duer A\duer C^*)$ and
$G\in\cinf{C^*}$. 

From the proof of \ref{theorem:bothduals} we have, in particular,
that
\begin{equation}
\label{eq:specials}
\begin{split}
[\xi^\sqcap,\,\eta^\sqcap] = [\xi,\eta]^\sqcap,\qquad
[\xi^\sqcap,\,\Bar{\phi}] = \Bar{\Lambda_\xi(\phi)},\qquad
[\Bar{\phi},\,\Bar{\psi}] = 0,\qquad
e(\xi^\sqcap)(\ell_c) = \ell_{Q_\xi(c)},\\
e(\xi^\sqcap)(f\circ q_{C^*}) = x(f)\circ q_{C^*},\qquad
e(\Bar{\phi})(\ell_c) = -\langle\phi,\,\da_A(c)\rangle\circ q_{C^*},
\qquad
e(\Bar{\phi})(f\circ q_{C^*}) = 0,
\end{split}
\end{equation}
where $(\xi, x)$ and $(\eta,y)$ are linear sections of $D\to A$, 
$\phi, \psi\in\Ga A^*$, $c\in\Ga C$ and $f\in\cinf{M}$. The 
operators $\Lambda_\xi\co\Ga A^*\to\Ga A^*$ and 
$Q_\xi\co\Ga C\to\Ga C$ are defined in the course of the proof. 

Notice that, since the relationship between $D$ and 
$D\duer A\duer C^*$ is reciprocal, there is a similar set of
equations for the bracket on $D\to A.$

For any vector bundle $V\to M$ with a linear Poisson structure $\pi$,
the Poisson anchor $\pi^{\#}\co T^*V\to TV$ is a morphism
of double vector bundles over $V$ and the anchor $a\co V^*\to TM$
of the dual Lie algebroid structure. Taking the dual of $\pi^{\#}$
over $V$ (that is, the usual dual) exchanges the side and core 
structures with the other's duals, so that $(\pi^{\#})^*$ has core 
morphism $a^*$. Since $(\pi^{\#})^* = -\pi^{\#}$ it follows that 
the core morphism of $\pi^{\#}$ is $-a^*$. See 
\cite[9.2.1]{Mackenzie:GT}. 

This argument applies to each vector bundle structure on $D\duer A$. 
Denote the Poisson anchor associated to $\pi_{\du A}$ by $\pi^{\#A}.$
Figure \ref{fig:cofd} shows the cotangent and tangent triple vector 
bundles for $D\duer A$, which are the domain and target for 
$\pi^{\#A}.$

\begin{figure}[h]
\begin{picture}(340,150)(0,30)
\put(0,150){\xymatrix@=5mm{
T^*(D\duer A) \ar[rr] \ar[dd] \ar[dr] 
& & D\duer A\duer C^* \ar[dd] \ar[dr] &\\
& D  \ar[rr] \ar[dd]         & & B \ar[dd]  \\
D \duer A   \ar[rr] \ar[dr]         & & C^* \ar[dr] &\\
& A  \ar[rr] & & M\\
}}
\put(100,30){(a)}
\put(250,150){\xymatrix@=6mm{
T(D\duer A) \ar[rr] \ar[dd] \ar[dr] 
& & TC^* \ar[dd] \ar[dr] &\\
& TA  \ar[rr] \ar[dd]         & & TM \ar[dd]  \\
D\duer A    \ar[rr] \ar[dr]         & & C^* \ar[dr] &\\
& A  \ar[rr] & & M\\
}}
\put(320,30){(b)}
\end{picture}
\caption{\ \label{fig:cofd}}
\end{figure}

Since $\pi_{\du A}$ is linear over $A$, 
$\pi^{\#A}\co T^*(D\duer A)\to T(D\duer A)$ is a morphism
of the left hand faces in Figure~\ref{fig:cofd}, and 
the map $D\to TA$ is the anchor $\dan_A$. 

Likewise, $\pi^{\#A}$ is a morphism of double vector bundles 
for the rear faces of Figure~\ref{fig:cofd} and the corner
map $D\duer A\duer C^*\to TC^*$ is the anchor $e$ for the Lie 
algebroid structure on $D\duer A\duer C^*$. 

It follows by a simple argument that $\pi^{\#A}$ is a morphism of 
triple vector bundles. The corner map $B\to TM$ is the anchor 
$b$ and the other four corner maps are identities. 

In effect, Theorem \ref{theorem:bothduals} shows that dualizing an 
\LAvb{} along its Lie algebroid structure gives rise to a 
\emph{Poisson double vector bundle}; that is, to a double vector 
bundle with a Poisson structure for which the Poisson anchor is a 
morphism of triple vector bundles. The converse may also be proved. 

\section{ABSTRACT DOUBLE LIE ALGEBROIDS}
\label{sect:adla}

We now turn to the general notion of double Lie algebroid. Before 
giving the definition, we need some preliminaries. 

Again consider a double vector bundle as in 
Figure~\ref{fig:intro}(a). We now assume that there are Lie 
algebroid structures on all four sides and that each pair of
parallel sides forms an \LAvb. 

Denote the four anchors by $\dan_A\co D \to TA,\
\dan_B\co D\to TB,$ and $b\co B\to TM,\
a\co A\to TM$. As usual we denote all four brackets by 
$[\ ,\ ]$; the notation for elements will make clear which 
structure we are using.

The anchors thus give morphisms of double vector bundles
\begin{gather*}
(\dan_A;\id,b;\id)\co (D;A,B;M)\to (TA;A,TM;M),\\
(\dan_B;a,\id;\id)\co (D;A,B;M)\to (TB;TM,B;M)
\end{gather*}
and so define morphisms of their cores; denote these by
$\da_A\co C\to A$ and $\da_B\co C\to B$.

From \S\ref{sect:pc} we know that $D\duer A\duer C^*$ is a Lie
algebroid over $C^*$. Likewise, the Lie algebroid structure on
$D\to B$ induces a Poisson structure on $D\duer B$ which is
linear over $C^*$ and therefore induces a Lie algebroid structure
on $D\duer B\duer C^*$ with base $C^*$. 

The bundles $D\duer A\duer C^*$ and $D\duer B\duer C^*$ are dual
over $C^*$ under the pairing:
$$
\newp{\mscr{X}}{\mscr{Y}} = 
\langle \mscr{X},\,Z_A^{-1}(\mscr{Y})\rangle =
\langle Z_B^{-1}(\mscr{X}),\,\mscr{Y}\rangle. 
$$

\begin{df}                                    
\label{df:doubla}
A \emph{double Lie algebroid} is a double vector bundle as in
Figure~{\rm \ref{fig:intro}(a)}, equipped with Lie algebroid 
structures on all four sides such that each pair of parallel
sides forms an \LAvb, and such that the Lie algebroid structures
on $D\duer A\duer C^*\to C^*$ and $D\duer B\duer C^*\to C^*$ form a
Lie bialgebroid. 
\end{df}

In practice we will use $Z_A$ or $Z_B$ to transfer the Lie
algebroid structure from, for example, $D\duer B\duer C^*$ 
to $D\duer A$. 

We recall the notion of Lie bialgebroid, which was introduced 
in \cite{MackenzieX:1994} as an abstraction of the infinitesimal 
structure of a Poisson groupoid \cite{Weinstein:1988}. The 
conventions used here follow \cite{Mackenzie:GT}. 

\begin{df}
Let $A$ be a Lie algebroid on base $M$ and suppose that $A^*$ is also
equipped with a structure of Lie algebroid. These structures 
constitute a \emph{Lie bialgebroid} if 
\begin{equation}
\label{eq:bialgd}
d_*[X,Y] = [d_*X,Y] + [X,d_*Y]
\end{equation}
for all $X,Y\in\Ga A$. 
\end{df}

In this definition $d_*$ is the coboundary operator for $A^*$ 
and $[~,~]$ is the Schouten bracket of multisections of $A$. See 
\cite{Mackenzie:GT}, for example. 

A Lie bialgebra is a Lie bialgebroid for which $M$ is a single 
point. For a Poisson manifold $M$, the cotangent Lie algebroid 
structure on $T^*M$
makes $(T^*M,\,TM)$ a Lie bialgebroid. The following material, which
is needed below, can be found in \cite{MackenzieX:1994} and 
\cite[Chapter~12]{Mackenzie:GT}. 

Suppose that $(A, A^*)$ is a Lie bialgebroid on base $M$ and 
denote the anchors by $a$ and $a_*$. Then we take the 
Poisson structure on $M$ to be $\pi^\#_M = a_*\circ a^*$. 
It follows that $a$ is a Poisson map 
(to the tangent lift Poisson structure on $TM$; see below) 
and $a_*$ is anti--Poisson.

The concept of Lie bialgebroid is self--dual: if $(A,A^*)$ is a Lie
bialgebroid, then $(A^*,A)$ is also (\cite[3.10]{MackenzieX:1994}, 
\cite{Kosmann-Schwarzbach:1995}, or see \cite{Mackenzie:GT}).
However $(A^*,A)$ will induce on $M$ the opposite of the Poisson 
structure induced by $(A,A^*)$. 

In this section we need the following \emph{morphism criterion} for 
a Lie bialgebroid (\cite[6.2]{MackenzieX:1994}, or 
see \cite[\S12.2]{Mackenzie:GT}).  

\begin{thm}
\label{theorem:6.2}
Let $A$ be a Lie algebroid on $M$ such that its dual vector 
bundle $A^{*}$ also has a Lie algebroid structure. Denote their 
anchors by $a$ and  $a_*$.
Then $(A, A^{*})$ is a Lie bialgebroid if and only if
\begin{equation}                         \label{main}
T^*(A^*)\buildrel{R}\over\longrightarrow T^*(A)
\buildrel\pi^\#_A\over\longrightarrow TA
\end{equation}
is a Lie algebroid morphism over $a_*$, where the domain
$T^{*}(A^{*})\to A^{*}$ is the cotangent Lie algebroid induced by the
Poisson structure on $A^{*}$, and the target $TA\to TM$ is the 
tangent prolongation of $A$.
\end{thm}

The tangent prolongation Lie 
algebroid\index{tangent prolongation (Lie algebroid)} structure 
on $TA\to TM$ is best thought of in terms of the dual structure. 
For any Poisson manifold $P$ there is a 
\emph{tangent lift}\index{tangent lift (Poisson)} 
Poisson structure on $TP$ for which the anchor is given by
\begin{equation}
\label{eq:tangentP}
\pi^\#_{TP} = J_P\circ T(\pi^\#_P)\circ\Theta_P^{-1}
\end{equation}
where $J_P\co T^2P\to T^2P$ is the canonical 
involution\index{canonical involution} for the manifold 
$P$ and $\Theta_P\co T(T^*P)\to T^*(TP)$ was named the 
\emph{Tulczyjew diffeomorphism}\index{Tulczyjew diffeomorphism} 
\emph{of} $P$ in \cite[p.395]{Mackenzie:GT} ($\Theta$ was
denoted $J'$ in \cite{MackenzieX:1994}). For any manifold $P$, 
either of the two formulas
\begin{equation}
\label{eq:Theta}
\Theta_P = R_{TM}\circ (\delta\nu)^\flat = J_M^*\circ I_{TM}
\end{equation}
may be taken as the definition of $\Theta_P$. Here 
$\delta\nu$ is the canonical symplectic structure on $T^*P$.

If $P = A^*$ is the dual
of a Lie algebroid then the tangent lift Poisson structure is 
linear with respect to the bundle $T(A^*)\to TM$ and therefore 
induces a Lie algebroid structure on $TA\to TM$ via the tangent 
pairing $\llangle~,~\rrangle.$ The anchor is $J_M\circ T(a).$ 
See \cite[5.6]{MackenzieX:1994} or 
\cite[10.3.14, 9.7.1]{Mackenzie:GT}. 

With these preliminaries established, consider a double Lie
algebroid $(D;A,B;M)$. We use $Z_A$ to transfer the Lie algebroid
structure from $D\duer B\duer C^*$ to $D\duer A$. 
From the morphism criterion, we know that 
$$
\pi^{\# A}\circ R\co T^*(D\duer A\duer C^*)\to T(D\duer A)
$$
is a morphism of Lie algebroids over $e\co D\duer A\duer C^*\to
TC^*$. As with any Lie bialgebroid, it is also a morphism of
double vector bundles, with the other side morphism being the
identity of $D\duer A$ and the core morphism being
$e^*\co T^*C^*\to D\duer A$. 

From \S\ref{sect:pc} we know that $e$ is a morphism of 
double vector bundles, namely of the two right--hand faces in 
Figure~\ref{fig:cofd}, over $C^*$ and $b\co B\to TM$. 

\begin{lem}
\label{lem:dualt}
The core of $e\co D\duer A\duer C^*\to TC^*$ is 
$-\da_A^*\co A^*\to C^*.$
\end{lem}

\pf
The dual of $\pi^{\# A}$ over $D\duer A$ is $-\pi^{\# A}$. 
In the dualization process, the two triple vector bundles in
Figure~\ref{fig:cofd} are interchanged, but the positioning of the
sides and cores is permuted. In particular, the cores of the 
right--hand faces and of the front faces are exchanged. The map 
of the front cores induced by $\pi^{\# A}$ is $\da_A$, by 
definition.
\pfend

The process of dualizing a triple structure is covered in 
\cite{Mackenzie:2005dts}. Lemma \ref{lem:dualt} can also be proved 
by using (\ref{eq:specials}). 

Dualizing $e$ over $C^*$ it follows that $e^*$ is a morphism of
double vector bundles over $C^*$ and $-\da_A$, with core morphism
$b^*\co T^*M\to B^*$. 

The Poisson structure induced on $C^*$ by the Lie bialgebroid 
$(D\duer A\duer C^*,\,D\duer A)$ has associated anchor
$$
\pi^{\#}_{C^*} = e_*\circ e^*\co T^*C^*\to TC,
$$
where $e_*\co D\duer A\to TC^*$ is the anchor of $D\duer A\to C^*$. 
Denoting the anchor of $D\duer B\duer C^*\to C^*$ by $f$, we have
$e_* = f\circ Z_A.$ The arguments applied above to $e$ apply
equally to $f$ and so $f\co D\duer B\duer C^*\to TC^*$ is a
morphism of double vector bundles over $C^*$ and $-a$, with core
morphism $-\da_B^*$. Thus $\pi^{\#}_{C^*}$ is a morphism of
double vector bundles and therefore induces on $C$ the structure 
of a Lie algebroid, for which the anchor is the side morphism
$(-a)\circ(-\da_A).$

The core morphism for $\pi^{\#}_{C^*}$ is 
$(-\da_B^*)\circ b^*.$ The skew--symmetry of a Poisson
anchor implies that its core morphism is the negative dual of its
non--identity side morphism, and so we have 
$a\circ\da_A = b\circ\da_B$. 

\begin{prop}
\label{prop:da}
The maps $\da_A\co C\to A$ and $\da_B\co C\to B$ are Lie algebroid 
morphisms.
\end{prop}

We need two preparatory lemmas. The first is from
\cite[10.3.6]{Mackenzie:GT} and is proved by the same method as the
second. 

\begin{lem}
\label{lem:proj1}
For any Lie algebroid $V$ on $M,$ the bundle projection
$T^*V^*\to V$ is a morphism of Lie algebroids over the bundle
projection $V^*\to M.$
\end{lem}

\begin{lem}
\label{lem:proj2}
The projections $D\duer A\duer C^*\to B$ and 
$D\duer B\duer C^*\to A$ are morphisms of Lie
algebroids over the bundle projection $C^*\to M$. 
\end{lem}

\pf
Denote the projections by $Q\co D\duer A\duer C^*\to C^*$ and 
$q\co C^*\to M$. The anchor condition 
$b\circ Q = T(q)\circ f$ follows from the fact that the Poisson
anchor for $D\duer A$, namely $T^*(D\duer A)\to T(D\duer A)$, is
not only a morphism of double vector bundles, but a morphism of
triple vector bundles. 

Since $q$ is a surjective submersion and $Q$ is 
fibrewise surjective, it is sufficient to check that if 
$\mscr{X}$ projects to $x$ and $\mscr{Y}$ projects to $y$, then
$[\mscr{X},\,\mscr{Y}]$ projects to $[x,y].$ By \ref{prop:gen}, 
it is sufficient to consider the cases of linear sections and
core sections. 

For two linear sections $(\xi,x)$ and $(\eta,y)$ of $D\to A$
we have from the proof of Theorem \ref{theorem:bothduals} that 
$[\xi,\eta]$ is linear over $[x,y]$ and 
$[\xi^\sqcap,\,\eta^\sqcap] = [\xi,\eta]^\sqcap.$ Likewise, 
for $\phi,\psi\in \Ga A^*$, the core sections project to
zero and $[\Bar{\phi},\,\Bar{\psi}] = 0$ does also. 

Lastly, for a linear section $(\xi,x)$ and $\phi\in\Ga A^*$, we
have from (\ref{eq:specials}) that $[\xi^\sqcap,\,\Bar{\phi}]$ is 
the core section  $\Bar{\Lambda_\xi(\phi)}$, and this projects to
zero. 

The result for $D\duer B\duer C^*\to A$ is proved in the same way. 
\pfend

\noindent
{\sc Proof of \ref{prop:da}:}
As noted above, $e^*\co T^*C^*\to D\duer A$ is a morphism of double
vector bundles over $C^*$ and $-\da_A\co C\to A$. As for any Lie 
bialgebroid, $e^*$ is also a morphism of Lie algebroids over $C^*$. 
Now the projection $T^*C^*\to C$ is a morphism of Lie algebroids,
and the projection $D\duer A\to A$ is an anti--morphism of Lie
algebroids (because $Z_A$ is $-\id$ on $A$). Both projections are
surjective submersions and so determine the Lie algebroid
structures on their targets. It follows that $-\da_A$ is an
anti--morphism of Lie algebroids. 
\pfend

The general process used here, of quotienting a known morphism
of Lie algebroids, is based on the following result. A 
\emph{fibration of Lie algebroids} is a morphism $\nu\co A\to A'$ 
over $p\co M\to M'$ such that $p$ is a surjective submersion and
$\nu$ is fibrewise surjective. 

\begin{prop}
\label{prop:descent}
Consider Lie algebroids $A,\ A'$ and $B$ on bases $M,\ M'$ and $N$, 
and let $\phi\co A\to A',\ \psi\co A'\to B$ and $\nu\co A\to A'$ 
be morphisms of vector bundles over $f,g$ and $p$, 
with $\phi = \psi\circ\nu$. 

If $\phi$ is a morphism of Lie algebroids and $\nu$ is a
fibration of Lie algebroids, then $\psi$ is a morphism
of Lie algebroids. 
\end{prop}

This follows from the general theory of quotients of Lie
algebroids \cite[\S4.4]{Mackenzie:GT}, but is easy to prove
directly. 

Observe that Proposition \ref{prop:da} cannot be proved
by mapping sections of $C$ to core sections of $D\to A$. In
the case of an \LAvb{} there is a simple relation (equation 
(\ref{eq:pbacks})), but this does not hold for a general 
double Lie algebroid. This may be seen by considering the iterated
tangent bundle $T^2M = T(TM)$ of any manifold; the bracket of two
vertical lift vector fields (which are the core sections in this
case) is always zero. 

The next result, that the anchors $\dan_A$ and $\dan_B$ are
morphisms of Lie algebroids with respect to both structures, was 
taken as part of the definition of a double Lie algebroid in 
\cite{Mackenzie:notions}. I am grateful to Ted Voronov for having 
shown \cite{Voronov:mtqm}, using super methods, that it is a
consequence of the other conditions in the definition. The proof
which follows uses standard methods. 

\begin{thm}
\label{thm:condII}
The anchor $\dan_A\co D\to TA$ is a Lie algebroid morphism
over $b\co B\to TM$, and the anchor $\dan_B\co D\to TB$ is a Lie 
algebroid morphism over $a\co A\to TM$. 
\end{thm}

\begin{figure}[h]
\begin{picture}(340,150)(0,30)
\put(100,150){\xymatrix@=5mm{
T^*(D\duer A\duer C^*)\ar@{-->}[rr]^{\pi^{\#A}\circ R}\ar[dd]
                         \ar[dr]^(0.7)\nu
& & T(D\duer A) \ar[dd] \ar[dr]^(0.7){T(\nu')} &\\
& D  \ar@{-->}[rr]^(0.33){\dan_A} \ar[dd]   & & TA \ar[dd]  \\
D \duer A\duer C^*   \ar@{-->}[rr]_e \ar[dr] & & TC^* \ar[dr] &\\
& B  \ar@{-->}[rr]_b & & TM\\
}}
\end{picture}
\caption{\ \label{fig:anchors}}
\end{figure}

\pf
Figure \ref{fig:anchors} shows part of the outline of the map 
$\pi^{\#A}\circ R\co T^*(D\duer A\duer C^*)\to TC^*$. The 
dashed lines
indicate that this should be read as a morphism of double structures,
from left to right, and emphasize that it is not a triple structure, 

Note first that every vertical arrow is the bundle projection of
a Lie algebroid. Next, $\pi^{\#A}\circ R$ is a morphism of Lie 
algebroids over $e$, by the morphism criterion. Next, $v$ and
$T(\nu')$ are fibrations of Lie algebroids. That $\nu'\co D\duer A
\to A$ is a morphism over $q_*\co C^*\to M$ follows from Lemma
\ref{lem:proj2}, and so its tangent prolongation is also. That
$\nu$ is a fibration follows from \ref{lem:proj1} by using 
$D\duer A\duer C^*\isom D\duer B$ and taking $V = D\to B$. 

The result now follows by Proposition \ref{prop:descent}. 
\pfend

In the sections which follow we consider three classes 
of examples of double Lie algebroids. 

\section{THE COTANGENT DOUBLE OF A LIE BIALGEBROID}
\label{sect:dlalba}

The notion of Lie bialgebroid is recalled in the previous section. 

Consider a Lie algebroid $A$ on $M$ together with a Lie algebroid
structure on the dual, not \emph{a priori} related to that on $A$. 
The structure on $A^*$ induces a Poisson structure on $A$, and this 
gives rise to a cotangent Lie algebroid $T^*A\to A$. Equally, the 
Lie algebroid structure on $A$ induces a Poisson structure on $A^*$ 
and this gives rise to a cotangent Lie algebroid $T^*A^*\to A^*$. 
We transfer this latter structure to $T^*A\to A^*$ via $R$. There 
are now Lie algebroid structures on each of the four sides of the 
double vector bundle $D = T^*A$ of Figure~\ref{fig:bialgd}(a).

\begin{thm}
\label{theorem:Manin}
Let $A$ be a Lie algebroid on $M$ such that its dual vector bundle 
$A^{*}$ also has a Lie algebroid structure. Then $(A, A^*)$ is a 
Lie bialgebroid if and only if $D = T^*A$, with the structures just 
described, is a double Lie algebroid.
\end{thm}

\pf
First assume merely that $A^*$ is a Lie algebroid. Consider 
Figure~\ref{fig:intro}(c) with the corresponding Lie algebroid 
structures on the horizontal sides. 

\begin{lem}
Figure~{\rm \ref{fig:intro}(c)} with the Lie algebroid structures 
just described, is an \LAvb. 
\end{lem}

\pf
It must be proved that the vector bundle operations are Lie 
algebroid morphisms. That the bundle projection $T^*A\to A^*$ is a 
Lie algebroid morphism follows from Lemma \ref{lem:proj1}. That
the addition map is a Lie algebroid morphism follows from the
fact that $\delta\ell_\phi\add{A}\delta\ell_\psi = 
\delta\ell_{\phi+\psi}$ for $\phi,\psi\in\Ga A^*$. The other
conditions are similar. 
\pfend

Now assume that $(A, A^*)$ is a Lie bialgebroid. The vertical 
structure on $D$ is the cotangent Lie algebroid structure for the 
Poisson structure on $A$. Dualizing over $A$ we have 
Figure~\ref{fig:bialgd}(b) and the Poisson structure on $TA$ is 
the tangent lift of the Poisson structure on $A$. By \S\ref{sect:pc} 
this Poisson structure is also linear over $TM$. The dual 
$D\duer A\duer TM$ is shown in Figure~\ref{fig:bialgd}(c). The 
map $I_A\co T(A^*)\to T^\sol A$ described in Example 
\ref{ex:tangentd} is an isomorphism of double vector bundles to 
$(T(A^*);A^*,TM;M)$ and by \cite[10.3.14]{Mackenzie:GT} it is an 
isomorphism of Lie algebroids over $TM$ to the tangent prolongation
of $A^*\to M$. 

\begin{figure}[h]
\begin{picture}(340,100)(-30,0)
\put(0,90){\xymatrix@=15mm{
T^*A \ar[r] \ar[d] & A^* \ar[d] \\
A \ar[r] & M 
}}
\put(30,0){(a)}
\put(150,90){\xymatrix@=15mm{
TA \ar[r] \ar[d] & TM \ar[d] \\
A \ar[r] & M 
}}
\put(190,0){(b)}
\put(300,90){\xymatrix@=15mm{
T^\sol A\ar[r] \ar[d] & TM \ar[d] \\
A^* \ar[r] & M 
}}
\put(340,0){(c)}
\end{picture}
\caption{\ \label{fig:bialgd}}
\end{figure}

Now consider the horizontal Lie algebroid structure in
Figure~\ref{fig:bialgd}(a). This was transported via $R$ from
the cotangent Lie algebroid of $A^*$. By the same argument as in
the vertical case, the induced Lie algebroid structure on
$D\duer A^*\duer TM\isom TA\to TM$ is isomorphic to the tangent
prolongation Lie algebroid of $A\to M$. 

We must now prove that~:

\begin{prop}
\label{prop:Tbialg}
For $(A,A^*)$ a Lie bialgebroid on $M$, and the structures just
described, the pair $(TA, T^\sol A)$ is a Lie bialgebroid 
on $TM$. 
\end{prop}

\pf
We use the morphism criterion \ref{theorem:6.2}. We must
prove that
\begin{equation}                                \label{eq:comp}
T^*(T^\sol A)\buildrel R_\sol\over\longrightarrow T^*(TA)
\buildrel\pi^\#_{TA}\over\longrightarrow T^2A
\end{equation}
is a morphism of Lie algebroids over the anchor 
$J_M\circ T(a_*)\circ I_A^{-1}\co T^\sol A\to T^2M$ of $T^\sol A$. 

Here $R_\sol$ is the canonical map $R$ for $TA\to TM$. The domain 
of (\ref{eq:comp}) is the cotangent Lie algebroid for the Poisson 
structure on $T^\sol A$. The target of (\ref{eq:comp}) is the 
iterated tangent prolongation of the Lie algebroid structure of $A$.

We need two lemmas. The first is easily proved (see 
\cite[9.6.4]{Mackenzie:GT}). The $\Theta$ is the 
Tulczyjew diffeomorphism defined in (\ref{eq:Theta}). 

\begin{lem}
\label{lem:964}
Let $M$ be any manifold and write $p\co TM\to M$ and 
$c\co T^*M\to M$ for the bundle projections. Then 
for $\mscr{X}\in T(T^* M)$ and $\xi\in T(TM)$ with 
$T(c)(\mscr{X}) = T(p)(\xi)$, 
$$ 
\llangle\mscr{X},\,\xi\rrangle = 
\langle\Theta(\mscr{X}),\,J(\xi)\rrangle. 
$$ 
\end{lem}

In the next lemma, the map
$I_A^\times\co T^*(T^\sol A)\to T^*(TA^*)$ is the diffeomorphism
of the cotangents induced by the diffeomorphism $I_A$. 

\begin{lem}
\label{lem:Rsol}
$R_\sol = 
\Theta_A\circ T(R_A)\circ\Theta_{A^*}^{-1}\circ I_A^\times.$
\end{lem}

\pf
We begin by calculating $R_\sol$ in terms of (\ref{eq:vue}). 
Take $F\in T^*(T^\sol A)$ and $\xi\in T^2A$ projecting to
the same point of $TA$. Then
\begin{equation}
\label{eq:Rsol}
\langle R_\sol(F),\,\xi\rangle_{TA} = 
\llangle \mscr{X},\,\xi\rrangle_{T^2M} - 
\langle F,\,\mscr{X}\rangle_{T^\sol A},
\end{equation}
where $\mscr{X}\in T(T^\sol A)$ has the requisite projections. 
Write $\Phi = \Theta_{A^*}^{-1}(I_A^\times(F))\in T(T^*A^*).$
Then 
$$
\langle R_\sol((I_A^\times)^{-1}(\Theta_{A^*}(\Phi))),\,
\xi\rangle_{TA} = 
\llangle \mscr{X},\,\xi\rrangle_{T^2M} - 
\langle (I_A^\times)^{-1}(\Theta_{A^*}(\Phi)),\,
\mscr{X}\rangle_{T^\sol A},
$$
Define $\Psi = T(I_A)^{-1}(\mscr{X})\in T^2(A^*)$. Then the first 
term on the RHS is 
$$
\llangle T(I_A)(\Psi),\,\xi\rrangle_{T^2M} =
\left.\frac{d}{dt}%
\langle I_A(\Psi_t),\,\xi_t\rangle_{TM}\right|_0 =
\left.\frac{d}{dt}%
\llangle\Psi_t,\,\xi_t\rrangle_{TM}\right|_0.
$$
where $\Psi_t$ and $\xi_t$ are curves in $TA^*$ and $TA$ with
tangents $\Psi$ and $\xi$. Denote this expression by
$\lllangle\Psi,\,\xi\rrrangle$; this is a pairing of 
$T^2(A^*)$ and $T^2(A)$ over $T^2M$. 

The second term on the RHS of (\ref{eq:Rsol}) is equal to 
$\langle \Theta_{A^*}(\Phi),\,
T(I_A^{-1})(\mscr{X})\rangle_{T^\sol A}$. By Lemma \ref{lem:964}, 
this is $\llangle\Phi,\,J_{A^*}(\Psi)\rrangle_{TA^*}.$

Now consider $\Theta_A\circ T(R)$. Using \ref{lem:964} again, 
$$
\langle\Theta_A(T(R)(\Phi)),\,\xi\rangle_{TA}
= \llangle T(R)(\Phi),\,J_A(\xi)\rrangle_{TA} = 
\left.\frac{d}{dt}\langle R(\Phi_t),\,\eta_t\rangle \right|_0,
$$
where $\Phi_t$ is a curve in $T^*A^*$ with tangent $\Phi$ and 
$\eta_t$ is a curve in $TA$ with tangent $J_A(\xi)$. 
Now
$$
\langle R(\Phi_t),\,\eta_t\rangle = 
\llangle W_t,\,\eta_t\rrangle_{TM} - \langle\Phi_t,\,W_t\rangle_{A^*}
$$
where $W_t\in TA^*$ is any element with the required outline. A 
careful checking of the outlines of the various elements shows that 
we can take $W_t$ to be an integral curve of $J_{A^*}(\Psi).$ It is
then clear that 
$$
\left.\frac{d}{dt}\langle\Phi_t,\,W_t\rangle_{A^*}\right|_0
= \llangle\Phi,\,J_{A^*}(\Psi)\rrangle_{TA^*}.
$$
The result now follows from the next lemma. 
\pfend

\begin{lem}
\label{lem:lllrrr}
Let $A\to M$ and $A^*\to M$ be any dual vector bundles and let
$\lllangle~,~\rrrangle$ denote the pairing of $T^2(A^*)$ and
$T^2(A)$ over $T^2M$ used in the proof above. Then, for 
$\Psi\in T^2(A^*)$ and $\xi\in T^2(A)$ projecting to the same
element of $T^2M$, 
$$
\lllangle J_{A^*}(\Psi),\,J_A(\xi)\rrrangle = 
\lllangle\Psi,\,\xi\rrrangle.
$$
\end{lem}

\pf
Write $\xi = \da_s\,\da_t\,a(0,0)$ for a smooth $a(s,t)\in A$. 
The notation means that $a$ is defined in a neighbourhood of the
origin, and that $\xi$ is the tangent vector at 0 to the curve 
$s\mapsto \da_t\,a(s,0) = \frac{\da a}{\da t}(s,0)$ in $TM$. 

The projection condition ensues that we can likewise write 
$\Psi = \da_s\,\da_t\,b(0.0)$ for a smooth patch 
$b(s,t)\in A^*$ with $q(a(s,t)) = q_*(b(s,t))$ for all $(s,t)$ 
near $(0,0)$. Now 
$$
\lllangle\Psi,\,\xi\rrrangle = 
\left.\frac{d}{ds}\llangle\da_t\,a(s,0),\,
\da_t\,b(s,0)\rrangle\right|_{s = 0} = 
\left.\frac{\da^2}{\da s\da t}\langle a(s,t),\,
b(s,t)\rangle\right|_{(0,0)}
$$
and the result follows from the equality of mixed partials. 
\pfend

We now return to the proof of \ref{prop:Tbialg}. From Lemma
\ref{lem:Rsol} and (\ref{eq:tangentP}) we have
$$
\pi^\#_{TA}\circ R_\sol = 
J_A\circ T(\pi^\#_A\circ R_A)\circ\Theta_{A^*}^{-1}\circ I_A^\times. 
$$
Since $(A,A^*)$ is a Lie bialgebroid, 
$\pi^\#_A\circ R_A\co T^*A^*\to TA$ is a morphism of
Lie algebroids over $a_*$, so $T(\pi^\#_A\circ R_A)$ is a morphism 
of the prolongation structures over $T(a_*)$. Also, since $I_A$ is 
a Poisson diffeomorphism, $I_A^\times$ is an isomorphism of the
cotangent Lie algebroids over $I_A^{-1}.$ 

For any Poisson manifold, $\Theta_P\co T(T^*P)\to T^*(TP)$
is an isomorphism of Lie algebroids over $TP$ from the tangent 
prolongation of the cotangent Lie algebroid structure on $T^*P$ to 
the cotangent Lie algebroid of the tangent Poisson structure 
\cite[10.3.13]{Mackenzie:GT}. We apply this to $P = A^*$.

Finally we apply Lemma \ref{lem:J2}(i). Altogether, we have that
$\pi^\#_{TA}\circ R_\sol$ is a Lie algebroid morphism and this 
completes the proof of Proposition \ref{prop:Tbialg}.
\pfend

\begin{lem}
\label{lem:J2}
{\bf (i)} 
Let $A$ be a Lie algebroid on $M$ and let $T^2A\to T^2M$ be the
iterated prolongation Lie algebroid. Then the canonical 
involution $J_A\co T^2A\to T^2A$ is a Lie algebroid automorphism 
over $J_M$. 

{\bf (ii)} For any Poisson manifold $P$, the canonical involution
$J_P$ is a Poisson diffeomorphism, with respect to the iterated
tangent lift Poisson structure. 
\end{lem}

\pf
{\bf (i)} It is sufficient to show that the dual map 
$J_A\duer T^2M$ is a Poisson
diffeomorphism. The dual of $T^2A\to T^2M$ is $T^\sol(TA)\to T^2M$
and this is isomorphic to $T(T^\sol A)\to T^2M$ under the map
$I$ for $TA\to TM$. Composing with $T(I_A)$ we have that 
$T^\sol(TA)\to T^2M$ is isomorphic to $T^2(A^*)\to T^2M$ and with
this identification, $J_A$ becomes $J_{A^*}\co T^2(A^*)\to T^2(A^*).$
The result now follows from (ii). 

{\bf (ii)} 
Write $J = J_P$ and $J_T = J_{TP}.$ 
We must show that $T(J)\circ\pi^\#_{T^2P} = 
\pi^\#_{T^2P}\circ J^\times.$ Using the definition of
$\pi^\#_{T^2P}$, we have to show that
$$
T(J)\circ J_T\circ T(J)\circ J_T\circ T(J)\circ T^2(\pi^\#)
\circ T(\Theta)^{-1}\circ\Theta_T^{-1} = 
T^2(\pi^\#)\circ T(\Theta)^{-1}\circ \Theta_T^{-1}\circ J^\times.
$$
Now $T(J)\circ J_T\circ T(J)\circ J_T\circ T(J) = J_T$. If the
coordinates in $\R^3$ are $t_1,t_2,t_3$ then $T(J)$ corresponds
to the permutation $(2~3)$ and $J_T$ to the permutation $(1~2)$. 
Using this and the naturality property $J_T\circ T^2(\pi^\#) = 
T^2(\pi^\#)\circ J_{T^*P}$, we are reduced to showing that, 
for any manifold $P$, 
\begin{equation}
J_P^\times\circ \Theta_{TP}\circ T(\Theta_P) = 
\Theta_{TP}\circ T(\Theta_P)\circ J_{T^*P}.
\end{equation}
Take $\xi\in T^2(T^*P)$ and $\zeta\in T^3P$ with $T^2(c)(\xi) = 
p_{Y^2}(\zeta).$ Then, using \ref{lem:964}, 
$$
\langle \Theta_{T}\circ T(\Theta)\circ J_{T^*P} (\xi),\,\zeta\rangle
= \llangle T(\Theta)\circ J_{T^*P} (\xi),\, J_T(\zeta)\rrangle
$$
and
$$
\langle J^\times\circ \Theta_{T}\circ T(\Theta)(\xi),\,\zeta\rangle
= \llangle T(\Theta)(\xi),\, J_T\circ T(J)(\zeta)\rrangle
$$
and equality follows as in Lemma \ref{lem:lllrrr}. 
\pfend

Returning to the proof of \ref{theorem:Manin}, suppose that $A$ 
is a Lie algebroid on $M$ and that $A^*$ has a Lie algebroid 
structure, not \emph{a priori} related to the structure on $A$. 
Consider $D = T^*A$ with the two cotangent Lie algebroid structures 
arising from the Poisson structures on $A^*$ and $A$, and suppose 
that these structures make $D$ a double Lie algebroid.

Then in particular, by Theorem \ref{thm:condII}
the anchor $T^*A\to TA^*$ of the horizontal structure must be 
a Lie algebroid morphism over $a\co A\to TM$ with respect to the 
vertical structures. The anchor is
$$
T^*A\buildrel{R_{A^*}}\over\longrightarrow T^*A^*
\buildrel\pi^\#_{A^*}\over\longrightarrow TA^*
$$
That this be a Lie algebroid morphism over $a$ is
precisely the dual form of \ref{theorem:6.2}.

This completes the proof of Theorem \ref{theorem:Manin}.
\pfend

Recall the Manin triple characterization of a Lie bialgebra, as given
in \cite{LuW:1990}: Given a Lie bialgebra 
$({\mathfrak g}, {\mathfrak g}^*)$ the vector space direct sum 
${\mathfrak d} = {\mathfrak g}\oplus{\mathfrak g}^*$
has a Lie algebra bracket defined in terms of the two coadjoint
representations; this is the classical Drinfel'd double of 
$(\gog,\gog^*)$. This bracket is invariant under the pairing
$\langle X + \phi,\, Y + \psi\rangle = \langle\psi,\,X\rangle + 
\langle\phi,\,Y\rangle$ and both ${\mathfrak g}$ and 
${\mathfrak g}^*$ are coisotropic subalgebras. Conversely, if a Lie
algebra ${\mathfrak d}$ is a vector space direct sum 
${\mathfrak g}\oplus{\mathfrak h}$, both of which are coisotropic 
with respect to an invariant pairing of ${\mathfrak d}$ with 
itself, then ${\mathfrak h}\isom{\mathfrak g}^*$ and 
$({\mathfrak g}, {\mathfrak h})$ is a Lie bialgebra, with 
${\mathfrak d}$ as the double.

This result provided a characterization of the notion of Lie 
bialgebra in terms of a single Lie algebra structure on 
${\mathfrak d}$, the conditions being expressed in terms of the 
simple notion of pairing. It also demonstrated, as a simple 
consequence, that the notion of Lie bialgebra is self--dual:
$({\mathfrak g}, {\mathfrak g}^*)$ is a Lie bialgebra if and only 
if $({\mathfrak g}^*, {\mathfrak g})$ is so.

For Lie bialgebroids, Kosmann--Schwarzbach 
\cite{Kosmann-Schwarzbach:1995} gave an elegant proof of
self--duality, in terms of her criterion that $(A,A^*)$ is a Lie 
bialgebroid if and only if the coboundary of one structure is a 
derivation of the Schouten bracket of the other. 

For a Lie bialgebra $(\gog, \gog^*)$, the structure of the 
cotangent double Lie algebroid $T^*\gog\isom\gog\times\gog^*$ 
equips $T^*\gog$ with Lie algebroid structures on bases $\gog$ 
and $\gog^*$. These are the action Lie algebroids defined
by the coadjoint actions of $\gog$ and $\gog^*$ on each other. 
We show in \S\ref{sect:mpvdla} that a double Lie 
algebroid $(D;A,B.M)$ for which the core is trivial --- as is the 
case for $D = T^*\gog$ --- has a third Lie algebroid structure, on 
base $M$. In the case of $T^*\gog$, $M$ is a point and this third
structure is the Drinfel'd double Lie algebra. 

Despite the differences, it seems legitimate to regard Theorem 
\ref{theorem:Manin} as a generalization of the Manin triple theorem, 
since it provides a criterion, in terms of a notion of double, 
for Lie algebroid structures on a pair of dual bundles to 
constitute a Lie bialgebroid. 

In \cite{LiuWX:1997}, Liu, Weinstein and Xu place a bracket
structure on $A\oplus A^*$, where $(A, A^*)$ is a Lie bialgebroid, 
which is a direct generalization of the Lie algebra structure of
the Drinfel'd double, and they prove a direct generalization
of the Manin triple theorem. This is not a Lie algebroid structure,
but a Courant algebroid, a vector bundle with a bracket structure 
in which the Jacobi
and Leibniz identities do not necessarily hold. Whereas the 
cotangent double has two simple structures which interact in a 
nontrivial way, the single bracket structure on $A\oplus A^*$ of 
\cite{LiuWX:1997} is less well--behaved but embodies the whole
structure. 

The notion of Courant algebroid has been developed in work on
Dirac structures and generalized complex geometry. On the other 
hand, the notion of double Lie algebroid abstracts iterated and 
second--order constructions in 
differential geometry (see particularly \S\ref{sect:dladlg}), and 
provides a general setting for duality phenomena. In particular, 
there is a clear global form of the notion of double Lie 
algebroid. 

Another important property of the classical Drinfel'd double is
that $\gog\bowtie\gog^*$ itself has a natural structure of Lie 
bialgebra. A double form of this result will be given in 
\cite{MackenzieV}. 

\section{MATCHED PAIRS AND VACANT DOUBLE LIE ALGEBROIDS}
\label{sect:mpvdla}

A brief history of matched pairs of Lie algebras was given 
in the Introduction. The corresponding concept of a matched pair of 
Lie groups \cite{LuW:1990}, \cite{Majid:1990} was extended to 
Lie groupoids in \cite{Mackenzie:1992}. In \cite{Mokri:1997}, Mokri 
differentiated the twisted automorphism equations of 
\cite{Mackenzie:1992} to obtain conditions on a pair of Lie 
algebroid representations, of $A$ on $B$ and of $B$ on $A$,
which ensure that the direct sum vector bundle $A\oplus B$ has a Lie
algebroid structure with $A$ and $B$ as subalgebroids. We quote the
following.

\begin{df}  {\bf \cite[4.2]{Mokri:1997}}   
\label{df:mokri}
Let $A$ and $B$ be Lie algebroids on base $M$, with anchors $a$ 
and $b$, and let $\rho\co A\to\D(B)$ and $\sigma\co B\to\D(A)$ 
be representations of $A$ on the vector bundle $B$ and of $B$ on 
the vector bundle $A$. Then $A$ and $B$ together with $\rho$ 
and $\sigma$ form a \emph{matched pair of Lie algebroids} if 
equations $(\ref{eq:first})$ from the Introduction hold. 
\end{df}

The Lie algebroid $\D(E)$ of derivations on $E$ 
is defined in \S\ref{sect:pc}. The proof of the following proceeds 
along the same lines as the proof for Lie algebras. 

\begin{prop} {\bf \cite[4.3]{Mokri:1997}}        
\label{prop:mokri}
Given a matched pair of Lie algebroids, there is a Lie algebroid 
structure $A\bowtie B$ on the direct sum vector bundle 
$A\oplus B$, with anchor $c(X\oplus Y) = a(X) + b(Y)$ and bracket
$$ 
[X_1\oplus Y_1, X_2\oplus Y_2] =
 \{[X_1, X_2] + \sigma_{Y_1}(X_2) - \sigma_{Y_2}(X_1)\}\oplus
    \{[Y_1, Y_2] + \rho_{X_1}(Y_2) - \rho_{X_2}(Y_1)\}.
$$ 
Conversely, if $A\oplus B$ has a Lie algebroid structure for which
$A\oplus 0$ and $0\oplus B$ are Lie subalgebroids, then $\rho$ 
and $\sigma$ defined by 
$[X\oplus 0, 0\oplus Y] = -\sigma_Y(X) \oplus \rho_X(Y)$
form a matched pair.
\end{prop}

We now show that matched pairs of Lie algebroids correspond 
precisely to double Lie algebroids for which the core is zero. 

\begin{df}
A double Lie algebroid $(D; A, B; M)$ is \emph{vacant} if the
core is the zero bundle. 
\end{df}

Equivalently, a double Lie algebroid is vacant if the 
combination of the two projections,
$(q^D_A,\, q^D_B)\co D\to A\times_M B$ is a diffeomorphism.
This is a direct analogue of the definition of vacancy for 
double Lie groupoids and for \LAgpds{}
\cite[2.11, 4.10]{Mackenzie:1992}.

Consider a vacant double Lie algebroid $(D; A, B; M)$. Note that 
$D\to A$ and $D\to B$ are the pullback bundles $q_A^!B$ and 
$q_B^!A$. The dual $D\duer A \to C^*$ is the Whitney sum 
$A\oplus B^*$. The duality between $D\duer A$ and
$D\duer B$ is 
\begin{equation}
\label{eq:pairing}
\langle X + \psi,\, \phi + Y\rangle =
          \langle\psi,\, Y\rangle - \langle \phi,\, X \rangle. 
\end{equation}

The horizontal bundle projection $q^D_B\co D\to B$ is a
morphism of Lie algebroids over $q_A\co A\to M$ and since it is a
fibrewise surjection, it defines an action of $B$ on $q_A$ as in
\cite[\S4.1]{Mackenzie:GT}. Namely, each section $Y$ of $B$ induces 
the pullback section $1\otimes Y$ of $q_A^!B$ and this induces a 
vector field $\eta(Y) = \dan_A(1\otimes Y)$ on $A$, where 
$\dan_A\co D\to TA$ is the anchor of the vertical 
structure. Since the anchor is a morphism of double vector bundles, 
$\eta(Y)$ is linear over the vector field $b(Y)$ on $M$, in the 
sense of \cite[\S3.4]{Mackenzie:GT}; that is, $\eta(Y)$ is a vector 
bundle morphism $A\to TA$ over $b(Y)\co M\to TM$. It follows 
that $\eta(Y)$ defines derivations
$\sigma^{(*)}_Y$ on $A^*$ and $\sigma_Y$ on $A$ by
\begin{equation}                          
\label{eq:sigma}
\eta(Y)(\ell_\phi) = \ell_{\sigma^{(*)}_Y(\phi)},\quad
\langle\phi, \sigma_Y(X)\rangle = b(Y)\langle\phi, X\rangle
- \langle\sigma^{(*)}_Y(\phi), X\rangle
\end{equation}
where $\phi\in\Ga A^*,\ X\in\Ga A,$ and $\ell_\phi$ denotes the 
function $A\to\R,\ X\mapsto\langle\phi(q_AX), X\rangle$. Since 
$q^D_A$ is a Lie algebroid morphism, it follows that $\sigma$ is 
a representation of $B$ on the vector bundle $A$, with
contragredient representation $\sigma^{(*)}$. We will use
both $\sigma$ and $\sigma^{(*)}$ in what follows. 

Likewise, $q^D_A$ is a morphism of Lie algebroids over $q_B$ and 
for each $X\in\Ga A$ we obtain a linear vector field 
$\xi(X)\in \mathcal{X}(B)$ over $a(X)$. We define derivations 
$\rho^{(*)}_X$ on $B^*$ and $\rho_X$ on $B$ by
\begin{equation}                             
\label{eq:rho}
\xi(X)(\ell_\psi) = \ell_{\rho^{(*)}_X(\psi)},\quad
\langle\psi, \rho_X(Y)\rangle = a(X)\langle\psi, Y\rangle
- \langle\rho^{(*)}_X(\psi), Y\rangle.
\end{equation}
Again, $\rho^{(*)}$ and $\rho$ are representations of $A$.

The two Lie algebroid structures on $D$ are action Lie algebroids 
determined by the actions $Y\mapsto\eta(Y)$ and $X\mapsto\xi(X)$
\cite[\S4.1]{Mackenzie:GT}. It therefore follows that the dual 
Poisson structures are semi--direct, but we need to prove this
directly. We first calculate the Lie algebroid structure on
$A^*\oplus B$ induced from that on $D\to A.$

\begin{lem}                                  
\label{lem:sdps}
The Lie algebroid structure on $D\duer A\duer C^* = A^*\oplus B$ 
induced from the double Lie algebroid structure has
anchor $e(\phi\oplus Y) = b(Y)$ and bracket
\begin{equation}                              
\label{eq:sdpsigma}
[\phi_1\oplus Y_1,\, \phi_2\oplus Y_2] =
   \{\sigma^{(*)}_{Y_1}(\phi_2) 
       - \sigma^{(*)}_{Y_2}(\phi_1)\}\oplus[Y_1, Y_2].
\end{equation}
The Lie algebroid structure on $D\duer A = A\oplus B^*$ induced 
from the double Lie algebroid structure has 
anchor $e_*(X\oplus\psi) = -a(X)$ and bracket
\begin{equation}                            
\label{eq:sdprho}
[X_1\oplus\psi_1,\, X_2\oplus\psi_2] =
   [X_2, X_1]\oplus\{\rho^{(*)}_{X_2}(\psi_1) 
      - \rho^{(*)}_{X_1}(\psi_2)\}.
\end{equation}
\end{lem}

\pf
We apply equations (\ref{eq:specials}). Firstly, a linear section
$(\xi, Y)$ of $D\to A$ is necessarily a pullback section
$\xi(X) = X\oplus Y(m)$ for $X\in A_m.$ From the definition 
(\ref{eq:xisqcap}) of $\xi^\sqcap$ we find that
$\xi^\sqcap\in\Ga(A^*\oplus B)$ is
$$
\xi^\sqcap(m) = 0_m^{A^*}\oplus Y(m). 
$$
It follows from the first equation in (\ref{eq:specials}) that
$[0\oplus Y_1,\,0\oplus Y_2)] = 0\oplus [Y_1,Y_2]$ for 
$Y_1,Y_2\in\Ga B.$

Next, $\phi\in\Ga A^*$ defines $\Bar{\phi}\in\Ga(A^*\oplus B)$
where $\Bar{\phi}(m) = \phi(m)\oplus 0^B_m.$ It follows from 
(\ref{eq:specials}) that $[\phi_1\oplus 0,\,\phi_2\oplus 0]
= 0\oplus 0.$

The second equation in (\ref{eq:specials}) is
$[\xi^\sqcap,\,\Bar{\phi}] = \Bar{\Lambda_\xi(\phi)}.$ Comparing
(\ref{eq:Lambda}) with (\ref{eq:sigma}), we have $\Lambda = 
\sigma^{(*)}.$ So $[0\oplus Y,\,\phi\oplus 0] = 
\sigma^{(*)}_Y(\phi)\oplus 0$ and this completes the proof of 
(\ref{eq:sdpsigma}). 

For (\ref{eq:sdprho}) apply (\ref{eq:specials}) to 
$D\duer B\duer C^* = A\oplus B^*$ and transport the structure to
$D\duer A = A\oplus B^*$ by $Z_A(X,\psi) = (-X,\psi).$
\pfend

Thus $A^*\oplus B$ is the semi--direct product $A^*\pds B$ 
of $B$ with the vector bundle $A^*$ with respect to $\sigma^{(*)}$. 
(This is the simple semi--direct product of two Lie algebroids on 
the same base, as in \cite[p.270]{Mackenzie:GT}.) We denote the
structure on $A\oplus B^*$ by $A^{op}\sdp B^*$. 

\begin{thm}                                     
\label{theorem:mp}
Let $(D;A,B;M)$ be a vacant double Lie algebroid. Then the two
Lie algebroid structures on $D$ are action Lie algebroids
corresponding to actions which define representations $\rho$, of $A$
on $B$, and $\sigma$, of $B$ on $A$, with respect to which $A$ and 
$B$ form a matched pair.

Conversely, let $A$ and $B$ be a matched pair of Lie algebroids with
respect to representations $\rho$ and $\sigma$. Then the action of 
$A$ on $q_B$ induced by $\rho$ and the action of $B$ on $q_A$ 
induced by $\sigma$ define Lie algebroid structures on 
$D = A\times_M B$ with respect to which
$(D;A,B;M)$ is a vacant double Lie algebroid.
\end{thm}

\pf
We first need to prove the three equations (\ref{eq:first}). 
We start by applying the bialgebroid condition to $E = A^*\pds B$
and $E^* = A^{op}\sdp B^*.$ We use the dual
form of (\ref{eq:bialgd}), which in this case is~:
\begin{equation}                         
\label{eq:full}
d^{E^*}[\phi_1 \oplus Y_1,\, \phi_2\oplus Y_2] =
   [d^{E^*}(\phi_1\oplus Y_1),\, \phi_2\oplus Y_2] +
   [\phi_1\oplus Y_1,\, d^{E*}(\phi_2\oplus Y_2)]
\end{equation}
for all $\phi_1\oplus Y_1, \phi_2\oplus Y_2\in \Ga E$. By replacing
$\phi_2\oplus Y_2$ by $f\phi_2\oplus fY_2$, for $f\in\cinf{M}$, and
expanding out, we find, as in \cite[p.449]{Mackenzie:GT}, that
\begin{equation}                         
\label{eq:onef}
d^{E^*}[\phi \oplus Y,\, f] = 
[d^{E^*}(\phi\oplus Y),\, f] + [\phi\oplus Y,\, d^{E*}(f)]. 
\end{equation}
Now $d^{E*}(f) = -d^Af\oplus 0$. So, for the last term, we have 
$[\phi\oplus Y,\, d^{E*}(f)] = -\sigma^{(*)}_Y(d^Af)\oplus 0,$ and
pairing this with any $X\oplus\psi\in\Ga E^*$. we have
$$
\langle[\phi\oplus Y,\, d^{E*}(f)],\,X\oplus\psi\rangle = 
a(\sigma_Y(X))(f) - b(Y)(a(X)(f)).
$$
On the LHS, $[\phi\oplus Y,\,f] = e(\phi\oplus Y)(f) = b(Y)(f)$ and
so, pairing with $X\oplus\psi$ gives
$$
\langle d^{E^*}(b(Y)(f)),\,X\oplus\psi\rangle = 
e_*(X\oplus\psi)(b(Y)(f)) = -a(X)(b(Y)(f)).
$$

Consider the first term on the RHS of (\ref{eq:onef}). For
any $\xi 
\in\Ga\extt{2}{E}$ we have $[\xi,f] = -\iota_{df}(\xi),$
\cite[12.1.4]{Mackenzie:GT}, so pairing with $X\oplus\psi$ gives
$$
\langle[\xi,f],\,X\oplus\psi\rangle = -\xi(df,\,X\oplus\psi).
$$

We need the following lemma, which is a straightforward calculation.

\begin{lem}
\label{lem:dE*}
$$
d^{E^*}(\phi\oplus Y)(X_1\oplus\psi_1,\, X_2\oplus\psi_2) =
   -(d^A\phi)(X_1,\, X_2) + \langle\psi_1,\, \rho_{X_2}(Y)\rangle
            - \langle\psi_2,\, \rho_{X_1}(Y)\rangle.
$$
\pfend
\end{lem}

For the first term on the RHS of (\ref{eq:onef}). we therefore have
\begin{multline*}
\langle[d^{E^*}(\phi\oplus Y),\,f],\,X\oplus\psi\rangle = 
-d^{E^*}(\phi\oplus Y)(d^Ef,\,X\oplus\psi) = \\
-\langle d^Bf,\,\rho_X(Y)\rangle = 
-b(\rho_X(Y))(f). 
\end{multline*}
using $d^Ef = 0\oplus d^Bf.$ Putting these together, 
(\ref{eq:onef}) becomes
$$
-a(X)(b(Y)(f)) = -b(\rho_X(Y))(f) + a(\sigma_Y(X))(f)
- b(Y)(a(X)(f)),
$$
which proves the third equation in (\ref{eq:first}). 

Now return to the bialgebroid equation (\ref{eq:full}). 
Set $\phi_1 = \phi_2 = 0$, so that we have 
\begin{equation}                         
\label{eq:320}
d^{E^*}[0 \oplus Y_1,\, 0\oplus Y_2] =
   [d^{E^*}(0\oplus Y_1),\, 0\oplus Y_2] +
   [0\oplus Y_1,\, d^{E*}(0\oplus Y_2)]
\end{equation}
and evaluate this at $(0\oplus \psi)\wedge(X\oplus 0).$ 
The LHS of (\ref{eq:320}) is easily seen to be
$\langle\psi,\, \rho_{X}[Y_2, Y_1]\rangle.$ 

On the RHS, consider the second term first. Writing the bracket
as the Lie derivative of the multisection $d^{E^*}(0\oplus Y_2)$ 
with respect to $0\oplus Y_1$, we have 
\begin{equation}
\label{eq:shift2}
\begin{split}
&\langle \ld_{0\oplus Y_1}(d^{E^*}(0\oplus Y_2)),\,
   (0\oplus\psi)\wedge(X\oplus 0)\rangle\\
& =  
\ld_{0\oplus Y_1}\langle d^{E^*}(0\oplus Y_2),\,
       (0\oplus\psi)\wedge(X\oplus 0)\rangle
 - \langle d^{E^*}(0\oplus Y_2),\, 
        \ld_{0\oplus Y_1}((0\oplus\psi)\wedge(X\oplus 0))\rangle
\end{split}
\end{equation}
In the first term, 
$d^{E^*}(0\oplus Y_2)(0\oplus\psi,X\oplus 0) = 
\langle\psi,\,\rho_X(Y_2)\rangle,$ and 
$e(0\oplus Y_1) = b(Y_1)$, so we have
$b(Y_1)\langle\psi,\, \rho_X(Y_2)\rangle.$ 

For the second term, note first that
$\ld_{0\oplus Y}(0\oplus\psi) = 0\oplus\ld_Y(\psi)$
and $\ld_{0\oplus Y}(X\oplus 0) = \sigma_Y(X)\oplus 0.$
Using these, the second term on the RHS of (\ref{eq:shift2}) is
$$
\langle\ld_{Y_1}(\psi),\,\rho_X(Y_2)\rangle + 
\langle\psi,\,\rho_{\sigma_{Y_1}(X)}(Y_2)\rangle.
$$
Calculating the first term on the RHS of (\ref{eq:320}) in the
same way, we have
\begin{multline*}
\langle\psi,\, \rho_{X}[Y_2, Y_1]\rangle = 
- b(Y_2)\langle\psi,\, \rho_X(Y_1)\rangle
+ \langle\ld_{Y_2}(\psi),\,\rho_X(Y_1)\rangle 
+ \langle\psi,\,\rho_{\sigma_{Y_2}(X)}(Y_1)\rangle\\
+ b(Y_1)\langle\psi,\, \rho_X(Y_2)\rangle
- \langle\ld_{Y_1}(\psi),\,\rho_X(Y_2)\rangle 
- \langle\psi,\,\rho_{\sigma_{Y_1}(X)}(Y_2)\rangle.
\end{multline*}
Finally, $b(Y_1)\langle\psi,\, \rho_X(Y_2)\rangle
- \langle\ld_{Y_1}(\psi),\,\rho_X(Y_2)\rangle = 
\langle\psi,\,[Y_1,\rho_X(Y_2)]\rangle$ and we have the first 
equation of (\ref{eq:first}). The second equation is proved in
the same way. 

Now consider the second part of \ref{theorem:mp}. It is 
straightforward to prove that $(A\times_M B;A,B;M)$, with the
action Lie algebroid structures on $A\times_M B$, is an 
\LAvb{} both horizontally and vertically. We must verify the
bialgebroid condition (\ref{eq:full}).

Firstly, (\ref{eq:320}) at $(0\oplus \psi)\wedge(X\oplus 0)$ can
be proved by reversing the argument above. The case of any 
$(X_1 \oplus \psi_1,\, X_2 \oplus \psi_2) = 
(X_1 \oplus 0,\, 0 \oplus \psi_2) 
+ (0 \oplus \psi_1,\, X_2 \oplus 0)$ follows from the skew--symmetry 
of all the terms. Thus (\ref{eq:320}) holds. 

For (\ref{eq:full}) with $Y_1 = Y_2 = 0$ and any arguments, it is
straightforward to see that all terms are identically zero. It
remains to consider (\ref{eq:full}) with $\phi_1 = 0,\ Y_2 = 0$.   
For the LHS we have
$$
d^{E^*}[0\oplus Y,\,\phi\oplus 0] = 
-d^A(\sigma^{(*)}_Y(\phi))\oplus 0. 
$$
Expanding out $-d^A(\sigma^{(*)}_Y(\phi))(X_1,\,X_2)$ we have
\begin{multline}
\label{eq:LHS}
-a(X_1)b(Y)\langle\phi,\,X_2\rangle
+a(X_1)\langle\phi,\,\sigma_Y(X_2)\rangle
+a(X_2)b(Y)\langle\phi,\,X_1\rangle\\
-a(X_2)\langle\phi,\,\sigma_Y(X_1)\rangle
+b(Y)\langle\phi,\,[X_1,\,X_2]\rangle
-\langle\phi,\,\sigma_Y([X_1,\,X_2])\rangle.
\end{multline}
On the RHS the second term is 
$[0\oplus Y,\,d^{E^*}(\phi\oplus 0)] 
= [0\oplus Y,\,-d^A\phi\oplus 0] = -\sigma^{(*)}_Y(d^A\phi)\oplus 0$,
using \ref{lem:dE*} and (\ref{eq:sdpsigma}), and extending 
$\sigma^{(*)}_Y$ to the exterior square. Expanding out 
$\langle-\sigma^{(*)}_Y(d^A\phi),\,X_1\wedge X_2\rangle$, we have
\begin{multline}
\label{eq:RHS2}
- b(Y)a(X_1)\langle\phi,\,X_2\rangle 
+ b(Y)a(X_2)\langle\phi,\,X_1\rangle  
+ b(Y)\langle\phi,\,[X_1,\,X_2]\rangle\\ 
+ a(\sigma_Y(X_1))\langle\phi,\,X_2\rangle
- a(X_2)\langle\phi,\,\sigma_Y(X_1)\rangle
- \langle\phi,\,[\sigma_Y(X_1),\,X_2]\rangle\\
+ a(X_1)\langle\phi,\,\sigma_Y(X_2)\rangle
- a(\sigma_Y(X_2))\langle\phi,\,X_1\rangle
- \langle\phi,\,[X_1,\,\sigma_Y(X_2)]\rangle.
\end{multline}
The first term on the RHS, evaluated at $F_1\wedge F_2\in
\Ga\extt{2}{E^*}$, is
\begin{equation}
\label{eq:RHS1}
\begin{split}
\langle[d^{E^*}(0\oplus Y),\,\phi\oplus 0],\,F_1\wedge F_2\rangle
= 
\langle-\ld_{\phi\oplus 0}(d^{E^*}(0\oplus Y)),\,F_1\wedge F_2\rangle
\\
= 
\langle d^{E^*}(0\oplus Y),\,
          \ld_{\phi\oplus 0}(F_1\wedge F_2)\rangle,
\end{split}
\end{equation}
since $e(\phi\oplus 0) = 0.$ 
Now $\ld_{\phi\oplus 0}(X\oplus\psi) = 0\oplus\theta_{\phi,X}$
where $\theta_{\phi,X}\in B^*$ is the element such that
$$
\langle Z,\,\theta_{\phi,X}\rangle 
= \langle\sigma^{(*)}_Z(\phi),\,X\rangle
$$
for all $Z\in B.$ Put $F_i = X_i\oplus\psi_i$ in (\ref{eq:RHS1})
and write $\theta_i = \theta_{\phi,X_i}.$ Then (\ref{eq:RHS1})
becomes
\begin{equation*}
\begin{split}
\langle d^{E^*}&(0\oplus Y),\,
(0\oplus\theta_1)\wedge(X_2\oplus\psi_2)\rangle
+ \langle d^{E^*}(0\oplus Y),\,
(X_1\oplus\psi_1)\wedge(0\oplus\theta_2)\rangle\phantom{XXXXX}\\
& = a(X_2)\langle Y,\,\theta_1\rangle
-\langle 0\oplus Y,\, 0\oplus\rho^{(*)}_{X_2}(\theta_1)\rangle
-a(X_1)\langle Y,\,\theta_2\rangle
+\langle0\oplus Y,\,0\oplus\rho_{X_1}^{(*)}(\theta_2)\rangle\\
& = b(\rho_{X_2}(Y))\langle\phi,\,X_1\rangle
-\langle\phi,\,\sigma_{\rho_{X_2}(Y)}(X_1)\rangle
-b(\rho_{X_1}(Y))\langle\phi,\,X_2\rangle
+\langle\phi,\,\sigma_{\rho_{X_1}(Y)}(X_2)\rangle.
\end{split}
\end{equation*}

Adding this to (\ref{eq:RHS2}) gives (\ref{eq:LHS}),
using the first and last equations of (\ref{eq:first}). 
The proof of Theorem~\ref{theorem:mp} is complete. 
\pfend

In the process we have also proved the following result.

\begin{cor}                                              
\label{cor:sdp}
Let $A$ and $B$ be Lie algebroids on the same base $M$, and let
$\rho\co A\to\D(B)$ and $\sigma\co B\to\D(A)$ be
representations on the vector bundles underlying $B$ and $A$. 

Then $A$, $B$, $\rho$, $\sigma$ form a matched pair if and only if 
$(A^*\pds B,\, A^{op}\sdp B^*)$, with the Lie algebroid structures 
described in {\rm \ref{lem:sdps}}, is a Lie bialgebroid.
\end{cor}

Thus the three matched pair equations in (\ref{eq:first}) are 
embodied in the bialgebroid equation (\ref{eq:full}).
Theorem \ref{theorem:mp} actually shows the equivalence of three 
formulations: a vacant double Lie algebroid structure with sides
$A$ and $B$ is equivalent to a matched pair structure on $(A, B)$, 
which in turn is equivalent to a Lie bialgebroid structure on 
$(A^*\pds B,\, A^{op}\sdp B^*)$ of the type in \ref{lem:sdps}. 
We thus have a diagrammatic characterization of matched pairs of 
Lie algebroids: recall that the 
characterization \ref{theorem:6.2} of a Lie bialgebroid is
formulated entirely in terms of the Poisson tensor and the canonical
isomorphism $R$. Thus \ref{theorem:mp} provides a definition of 
matched pair which can be formulated more generally in a category 
possessing pullbacks and suitable additive structure. The
corresponding characterization of matched pairs of Lie groupoids
in terms of vacant double Lie groupoids is in 
\cite[\S2]{Mackenzie:1992}. 

Given a vacant double Lie algebroid $(D;A,B;M)$ the Lie algebroid
structure $A\bowtie B$ on $D\to M$ may be called a \emph{diagonal
structure} by analogy with the groupoid case. 

Forms of Corollary \ref{cor:sdp} in the Lie algebra case have 
been found by Majid \cite[\S8.3]{Majid}
and in the setting of Lie--Rinehart algebras by 
Huebsch\-mann \cite{Huebschmann:2000}. 

In the case of a Lie bialgebra $({\mathfrak g}, {\mathfrak g}^*)$, 
the Lie bialgebroid associated to the vacant double Lie algebroid is
$(\gog\oplus\gog_0, \gog^*\oplus\gog^*_0)$ where the subscripts 
denote the abelianizations. This is of course consistent
with \ref{theorem:Manin} in the bialgebra case --- which is both a
bialgebroid and a matched pair.

\section{THE DOUBLE LIE ALGEBROID OF A DOUBLE LIE GROUPOID}
\label{sect:dladlg}

Definition \ref{df:doubla} of a double Lie algebroid
arose from work on the double Lie algebroid of a double Lie 
groupoid, as constructed in \cite{Mackenzie:1992}, 
\cite{Mackenzie:2000}. Some of what is required in order to verify 
that the double Lie algebroid of a double Lie groupoid does satisfy 
\ref{df:doubla} has been given in \cite{Mackenzie:1999}, and we 
recall the details briefly.

First we recall the notion of double Lie 
groupoid (see \cite{Mackenzie:1992} and references 
given there). A double Lie groupoid consists of a manifold $S$ 
equipped with two Lie groupoid structures on bases $H$ and $V$, 
each of which is a Lie groupoid on base $M$, such that 
the structure maps (source, target, multiplication, identity, 
inversion) of each groupoid structure on $S$ are morphisms with 
respect to the other, and such that the \emph{double source map}
$\alpha_2\co S\to H\times_M V,\ s\mapsto (\tilalpha_V(s),\,
\tilalpha_H(s))$ is a surjective submersion; see 
Figure~\ref{fig:S}(a).
\begin{figure}[h]
\begin{picture}(340,100)(-30,0)
\put(0,90){\xymatrix@=1.5cm{
S \ar@<0.5ex>[r]^{\tilalpha_H,\,\tilbeta_H}\ar@<-0.5ex>[r]
  \ar@<-0.5ex>[d]_{\tilalpha_V,\,\tilbeta_V}\ar@<0.5ex>[d]
& V \ar@<-0.5ex>[d]\ar@<0.5ex>[d]^{\alpha_V,\,\beta_V}\\
H \ar@<-0.5ex>[r]_{\alpha_H,\,\beta_H}\ar@<0.5ex>[r] & M\\
}}
\put(60,0){(a)}
\put(150,90){\xymatrix@=1.5cm{
A_VS \ar@<0.5ex>[r]^{A(\tilalpha_H),\,A(\tilbeta_H)}\ar@<-0.5ex>[r]
  \ar[d]_{\tilq_V}
& AV \ar[d]^{q_V}\\
H \ar@<-0.5ex>[r]\ar@<0.5ex>[r] & M\\
}}
\put(190,0){(b)}
\put(300,90){\xymatrix@=1.5cm{
A^2S \ar[r]^{\dot{q}_H}\ar[d]_{A(\tilq_V)}
& AV \ar[d]\\
AH \ar[r]_{q_H} & M\\
}}
\put(350,0){(c)}
\end{picture}
\caption{\ \label{fig:S}}
\end{figure}
One should think of elements of $S$ as squares, the horizontal 
edges of which come from $H$, the vertical edges from $V$, and 
the corner points from $M$.

Consider a double Lie groupoid $(S;H,V;M)$ as in 
Figure~\ref{fig:S}(a). Applying the Lie functor to the vertical 
structure $S\gpd H$ gives a Lie algebroid $A_VS\to H$ which also has 
a groupoid structure over $AV$ obtained by applying the Lie functor 
to the structure maps of $S\gpd V$; this is
the \emph{vertical \LAgpd}\ of $S$ \cite[\S4]{Mackenzie:1992}, as in
Figure~\ref{fig:S}(b). We recall the notion of \LAgpd{} from 
\cite{Mackenzie:1992}.

\begin{df}
An \emph{\LAgpd} is a manifold $\Om$ equipped with a Lie groupoid
structure on a base $A$, which is a Lie algebroid on a base $M$, 
and a Lie algebroid structure on a base $G$, which is a Lie
groupoid on $M$, such that the groupoid structure maps 
(source, target, identity, division) are morphisms of Lie algebroids,
and the map $\Om\to G\times_M A$ which combines the bundle projection
and the Lie groupoid source, is a surjective submersion. 
\end{df}

The Lie algebroid of $A_VS\gpd AV$ is denoted 
$A^2S$; there is a double vector bundle structure $(A^2S;AH,AV;M)$ 
obtained by applying $A$ to the vector bundle structure of 
$A_VS\to H$ \cite{Mackenzie:2000}; see Figure~\ref{fig:S}(c). 
Reversing the order of these operations, one defines first the 
\emph{horizontal \LAgpd}\ $(A_HS;AH,V;M)$ and then takes the Lie 
algebroid $A_2S = A(A_HS)$. The canonical involution 
$J_S\co T^2S\to T^2S$ then restricts to an isomorphism of double 
vector bundles $\jd_S\co A^2S\to A_2S$ and allows the Lie algebroid 
structure on $A^2S\to AV$ to be transported to $A_2S\to AV$. Thus 
$A_2S$ is a double vector bundle equipped with four Lie algebroid 
structures; in \cite{Mackenzie:2000} we called this the \emph{double 
Lie algebroid of $S$}. The core of both double vector bundles
$A_2S$ and $A^2S$ is $AK$, the Lie algebroid of the core groupoid 
$K\gpd M$ of $S$ \cite[1.6]{Mackenzie:2000}.

\begin{thm}
\label{thm:dladlg}
With the structures just defined, $(A_2S;AH,AV;M)$ is a double
Lie algebroid.   
\end{thm}

\pf
Write $D = A_2S$. The structure maps for the horizontal vector 
bundle $A_2S\to AV$ are obtained by applying the Lie functor to the 
structure maps of $A_HS\to V$ and are therefore Lie algebroid 
morphisms with respect to the vertical Lie algebroid structures. 
The corresponding statement is true for the vertical vector bundle 
$A^2S\to AH$ and this is transported by $\jd_S$ to $A_2S\to AH$. 
Thus the first condition of \ref{df:doubla} is satisfied. 

Now consider the bialgebroid condition. We first need to identify the
Lie algebroid structures on $A_2S\duer AH\duer A^*K$ and
$A_2S\duer AV\duer A^*K$. From \cite{Mackenzie:1999} we know 
that $A_V^*S\gpd A^*K$ is a Poisson groupoid, and therefore 
$(A(A_V^*S),\,A^*(A_V^*S))$ is a Lie bialgebroid. To relate these
structures, we need to recall the duality of \VBgpds. 

Consider any \VBgpd{} as in Figure \ref{fig:LAgpd}(a), with
core vector bundle $L\to M$. Taking the dual in the sense of 
\VBgpds{} \cite{Pradines:1988}, \cite[\S11.2]{Mackenzie:GT}, 
we have the \VBgpd{} $(\Om^*;G,L^*;M)$ as in 
Figure~\ref{fig:LAgpd}(b). 

\begin{figure}[htb]
\begin{picture}(340,100)(-30,0)
\put(0,90){\xymatrix@=1.5cm{
\Om \ar@<0.5ex>[r]\ar@<-0.5ex>[r] \ar[d] & A \ar[d]\\
G \ar@<-0.5ex>[r]\ar@<0.5ex>[r] & M\\
}}
\put(30,0){(a)}
\put(150,90){\xymatrix@=1.5cm{
\Om^* \ar@<0.5ex>[r]\ar@<-0.5ex>[r] \ar[d]
& L^* \ar[d]\\
G \ar@<0.5ex>[r]\ar@<-0.5ex>[r] & M\\
}}
\put(180,0){(b)}
\put(300,90){\xymatrix@=1.5cm{
A(\Om^*) \ar[r]\ar[d]
& L^* \ar[d]\\
AG \ar[r] & M\\
}}
\put(340,0){(c)}
\end{picture}
\caption{\ \label{fig:LAgpd}}
\end{figure}

The pairing of $\Om$ and $\Om^*$ is a Lie groupoid morphism
$\Om\times_G\Om^*\to\R$ and the induced pairing 
$\llangle~,~\rrangle = A(\langle~,~\rangle)$ of the Lie algebroids
is also non--degenerate and respects the double vector bundle 
structures \cite[\S11.5]{Mackenzie:GT}. It therefore 
induces an isomorphism of double vector bundles 
$$
\Goi_\Om\co A(\Om^*) \to A\Om\duer AG
$$
which preserves both side bundles and the core. Write $\Tilde{I}_V$
and $\Tilde{I}_H$ for this map when $\Om = A_VS$ and $\Om = A_HS.$
Composing the inverse of $\jd_S\duer AH$ with 
$\Tilde{I}_V$ gives an isomorphism
\begin{equation}
\label{eq:tdH}
\td_H = 
(\jd_S\duer AH)^{-1}\circ\Tilde{I}_V\co A(A_V^*S)\to 
A_2S\duer AH. 
\end{equation}
Dualizing this over $A^*K$ we get 
$\td_H\duer A^*K \co A_2S\duer AH\duer A^*K \to A^*(A_V^*S)$ and
to prove that this is an isomorphism of Lie algebroids it it
sufficient to prove that $\td_H$ is a Poisson map. 

On the other hand, consider $\Tilde{I}_H\co A(A_H^*S)\to
A_2S\duer AV$. Dualizing this over $A^*K$ we obtain
$$
\Tilde{I}_H\duer A^*K\co A_2S\duer AV\duer A^*K\to A^*(A_H^*S).
$$
Again, to show that this is an isomorphism of Lie algebroids, it
is sufficient to show that $\Tilde{I}_H$ is a Poisson map. However
the Poisson structure on $A_2S\duer AV$ is the dual of the Lie
algebroid structure on $A_2S\to AH$ and this is obtained from the
Lie algebroid structure on $A^2S\to AH$, transported by $\jd_S.$ 
So we actually need to show that $\td_V$, defined by 
\begin{equation}
\td_V = 
(\jd_S\duer AV)\circ\Tilde{I}_H\co A(A_H^*S)\to 
A^2S\duer AV
\end{equation}
is a Poisson map. The proof is therefore reduced to the
following. 

\begin{prop}
\label{prop:Theta}
The map $\td_H\co A(A_V^*S)\to A_2S\duer AH$ is a Poisson map with
respect to the Poisson structure on $A(A_V^*S)$ induced from the 
Poisson structure on $A_V^*S$ which is dual
to the Lie algebroid $A_VS\to H$, and the Poisson structure on
$A_2S\duer AH$ dual to the Lie algebroid $A_2S\to AH$. 

The map $\td_V\co A(A_H^*S)\to A^2S\duer AV$ is a Poisson map with 
respect to the Poisson structure on $A(A_H^*S)$ induced from the 
Poisson structure on $A_H^*S$ which is dual
to the Lie algebroid $A_HS\to V$, and the Poisson structure on
$A^2S\duer AV$ dual to the Lie algebroid $A^2S\to AV.$
\end{prop}

\pf
We will deduce these from the commutativity of 
Figures \ref{fig:IRTh}(a) and (b). 

\begin{figure}[h]
\begin{picture}(340,100)(-30,0)
\put(0,90){\xymatrix@=1.5cm{
A^*(A_H^*S) \ar[r]^{\Tilde{R}_H} \ar[d]_{\Lambda_H} 
& A_2S\duer AH \ar[d]^{\jd_S\duer AH}\\
A(A_V^*S) \ar[r]_{\Tilde{I}_V} \ar[ru]_{\td_H}
& A^2S\duer AH
}}
\put(60,0){(a)}
\put(200,90){\xymatrix@=1.5cm{
A^*(A_V^*S) \ar[r]^{\Tilde{R}_V} \ar[d]_{\Lambda_V} 
& A^2S\duer AV \\
A(A_H^*S) \ar[r]_{\Tilde{I}_H} \ar[ru]_{\td_V}
& A_2S\duer AV \ar[u]_{\jd_S\duer AV}
}}
\put(260,0){(b)}
\end{picture}
\caption{\ \label{fig:IRTh}}
\end{figure}

In turn we will deduce these diagrams from the commutative diagram
in Figure \ref{fig:IRTh-ori}(a), which is valid for any manifold $S$
\cite[p.442]{Mackenzie:GT}. Here $(\delta\nu)^\#$ is the map from
forms to fields induced by the canonical symplectic structure
$\delta\nu$ on $T^*S$. The lower right triangle is the definition
of $\Theta_S$ and the commutativity of the upper triangle is 
proved in \cite[9.6.7]{Mackenzie:GT}. 

\begin{figure}[h]
\begin{picture}(340,100)(-30,0)
\put(0,90){\xymatrix@=1.5cm{
T^*T^*S \ar[r]^{R_{TS}} \ar[d]_{(\delta\nu)^\#}
& T^*TS \ar[d]^{J_S^*}\\
TT^*S \ar[r]_{I_{TS}} \ar[ru]_{\Theta_S}
& T^\sol TS
}}
\put(60,0){(a)}
\put(200,90){\xymatrix@=1.5cm{
T^*T^*S \ar[r]^{(\delta\nu)^\#}
  \ar@<-0.5ex>[d]\ar@<0.5ex>[d]
& TT^*S \ar@<-0.5ex>[d]\ar@<0.5ex>[d]\\
A_V^*T^*S \ar[r]_{a_*} & 
TA_V^*S\\
}}
\put(260,0){(b)}
\end{picture}
\caption{\ \label{fig:IRTh-ori}}
\end{figure}

For any \VBgpd{} $(\Om;G,A;M)$ there is a canonical isomorphism of
double vector bundles $\goR_\Om\co A^*\Om^*\to A^*\Om$ which 
preserves the side bundles $A$ and $L^*$ and induces $-\id$ on the 
cores $A^*G$ \cite[11.5.13]{Mackenzie:GT}. The maps $\tilR_V$ and 
$\tilR_H$ in Figure \ref{fig:IRTh} 
are $\goR_\Om$ for $\Om = A^*_VS$ and $\Om = A^*_HS$. 

The maps $\Lambda_H$ and $\Lambda_V$ are the maps denoted 
$\mathcal{D}_H$ and $\mathcal{D}_V$ in \cite[\S3]{Mackenzie:1999}. We
recall the definition. 

The cotangent groupoid of the groupoid $S\gpd H$ is 
$T^*S\gpd A_V^*S$. This is a symplectic groupoid and so the induced 
map $(\delta\nu)^\#\co T^*T^*S\to TT^*S$ is an isomorphism of 
Lie groupoids as shown in Figure~\ref{fig:IRTh-ori}(b), where $a_*$ 
denotes the anchor for the Lie algebroid structure on 
$A_V^*T^*S\to A_V^*S.$

\begin{figure}[h]
\begin{picture}(340,100)(-30,0)
\put(0,90){\xymatrix@=1.5cm{
A_V^*T^*S \ar@<0.5ex>[r]\ar@<-0.5ex>[r]\ar[d]
& A^*T^*K \ar[d]\\
A_V^*S \ar@<-0.5ex>[r]\ar@<0.5ex>[r] & A^*K\\
}}
\put(60,0){(a)}
\put(150,90){\xymatrix@=1.5cm{
TA_V^*S \ar@<0.5ex>[r]\ar@<-0.5ex>[r]
\ar[d] & TA^*K \ar[d]\\
A_V^*S \ar@<0.5ex>[r]\ar@<-0.5ex>[r] & A^*K\\
}}
\put(200,0){(b)}
\put(300,90){\xymatrix@=1.5cm{
A_V^*TS \ar@<0.5ex>[r]\ar@<-0.5ex>[r]
\ar[d] & A^*TK \ar[d]\\
A_V^*S\ar@<0.5ex>[r]\ar@<-0.5ex>[r] & A^*K\\
}}
\put(350,0){(c)}
\end{picture}
\caption{\ \label{fig:dnu}}
\end{figure}

This $a_*$ is itself a morphism of double structures (indeed, of
triple structures); the domain and target are shown in 
Figures \ref{fig:dnu}(a) and (b). The \VBgpd{} (a) arises by
taking the vertical Lie algebroid of the double Lie groupoid
$(T^*S;A_V^*S,A_H^*S;A^*K)$ and then dualizing; it has core
$A^*(A_H^*S)$. The \VBgpd{} (b) is the tangent of $A_V^*S\gpd A^*K$
and so has core $A(A_V^*S)$. We define $\Lambda_H$ to be the core
map of $a_*.$ The definition of $\Lambda_V$ is similar. 

Both $R_{TS}$ and $\Theta_S$ are also groupoid morphisms. For any
\VBgpd{} $(\Om;G,A;M)$, $R_\Om$ is a groupoid morphism over
$\goR_\Om\co A^*\Om^*\to A^*\Om$ \cite[11.5.14]{Mackenzie:GT}. In
particular, $R_{TS}$ is a Lie groupoid morphism over 
$\goR_{TS/TH}\co A_V^*T^*S\to A_V^*TS$ where the subscript 
indicates that $\goR$ refers to the \VBgpd{} $(TS;TH,S;H).$

Similarly, $\Theta_S$ is a groupoid morphism over 
$\vartheta_{S/H}\co TA_V^*S\to A_V^*TS$ 
\cite[11.5.11]{Mackenzie:GT}. From the commutativity of
Figure \ref{fig:IRTh-ori}(a), it follows that 
$$
\vartheta_{S/H}\circ a_* = \goR_{TS/TH}.
$$
It further follows that the core maps of
$\vartheta_{S/H}$, $a_*$ and $\goR_{TS/TH}$ commute. The core
map of $a_*$ is $\Lambda_H$ by definition. To see the core of
Figure \ref{fig:dnu}(c), consider the triple structure in
Figure \ref{fig:triple}(a). The core of the top face here is
$A_VA_HS = A_2S$ and the core of the bottom face is $AH$. Taking the
dual of $A_VTS$ over $TH$ produces the triple structure shown in 
Figure \ref{fig:triple}(b). The bundle $A_2S\to AH$ of the top and
bottom cores of (a) is parallel to the axis of dualization and so
the bundle of top and bottom cores of (b) is $A_2S\duer AH\to AH$. 
The process of dualizing a triple vector bundle is described in 
detail in \cite{Mackenzie:2005dts} and the modifications to take 
account of groupoid structures are straightforward. 

\begin{figure}[h]
\begin{picture}(340,150)(0,30)
\put(0,150){\xymatrix@=5mm{
A_VTS \ar@<0.5ex>[rr]\ar@<-0.5ex>[rr] \ar[dd] \ar[dr] 
& & ATV\ar[dd] \ar[dr] &\\
& A_VS  \ar@<0.5ex>[rr]\ar@<-0.5ex>[rr] \ar[dd] & & AV \ar[dd]  \\
TH  \ar@<0.5ex>[rr]\ar@<-0.5ex>[rr] \ar[dr] & & TM \ar[dr] &\\
& H  \ar@<0.5ex>[rr]\ar@<-0.5ex>[rr] & & M\\
}}
\put(100,30){(a)}
\put(250,150){\xymatrix@=6mm{
A_V^*TS \ar@<0.5ex>[rr]\ar@<-0.5ex>[rr] \ar[dd] \ar[dr] 
& & A^*TK\ar[dd] \ar[dr] &\\
& A_V^*S \ar@<0.5ex>[rr]\ar@<-0.5ex>[rr] \ar[dd] & & A^*K \ar[dd]  \\
TH  \ar@<0.5ex>[rr]\ar@<-0.5ex>[rr] \ar[dr] & & TM \ar[dr] &\\
& H  \ar@<0.5ex>[rr]\ar@<-0.5ex>[rr] & & M\\
}}
\put(320,30){(b)}
\end{picture}
\caption{\ \label{fig:triple}}
\end{figure}

It is routine to verify that the core map of 
$\goR_{TS/TH}$ is $\tilR_H$ and the core map of $\vartheta_{S/H}$ 
is $\Tilde{\Theta}_H$. This proves the commutativity of the upper
triangle in Figure \ref{fig:IRTh}(a). The proof for (b) is similar. 

Now $R_{TS}\co T^*T^*S\to T^*TS$ is an antisymplectomorphism and
so $\goR_{TS/TH}$ is anti--Poisson \cite[11.5.16]{Mackenzie:GT}. 
Likewise $(\delta\nu)^\#\co T^*T^*S\to TT^*S$ is an 
antisymplectomorphism and so $\vartheta_{S/H}$ is anti--Poisson. 
It follows that the core maps are also anti--Poisson, and so 
$\Theta_H$ is a Poisson map. The case of $\Theta_V$ follows in the
same way. This concludes the proof of Proposition \ref{prop:Theta} 
and with it the proof of Theorem \ref{thm:dladlg}. 
\pfend

It is also possible to prove that for any \LAgpd{} $(\Om;G,A;M)$
applying the Lie functor yields a double Lie algebroid
$(A\Om;AG,A;M)$. In this case there is no version of the canonical
involution $J$ and it is necessary to work with prolonged Lie
algebroid structures; see \cite[3.3]{Mackenzie:2000}. 

Several examples of the double Lie algebroid of specific double
Lie groupoids are given in \cite[\S4]{Mackenzie:2000}. 
We conclude by noting three other examples. 

\begin{ex}\rm
Take $S$ to be the double groupoid $(M^4;M\times M,M\times M;M)$ in
which a quadruple of points is regarded as defining a square by its
four vertices. Then $A_2S$ is the double tangent bundle $T^2M$ and 
$\Tilde{\Theta}_H$ reduces to $\Theta_M.$ The associated Lie
bialgebroid is $(T^*T^*M,\,TT^*M)$ on base $T^*M$ with its
standard symplectic structure. 
\end{ex}

\begin{ex}\rm
An \emph{affinoid} $W$ \cite{Weinstein:1990} on a manifold $M$ 
may be regarded as a vacant double Lie subgroupoid of $M^4$ in the
previous example, such that both side groupoids $H$ and $V$ are 
the graphs of simple foliations defined by surjective submersions 
$\pi_1\co M\to Q_1$ and $\pi_2\co M\to Q_2$ 
\cite[\S3]{Mackenzie:1992}. The corresponding 
double Lie algebroid was calculated directly in 
\cite{Mackenzie:2000w} to be a double Lie subalgebroid
$(A_2W;T^{\pi_1}M,T^{\pi_2}M;M)$ of $T^2M$ with brackets defined by 
a pair of conjugate flat partial connections on $M$ adapted to the 
two foliations. The Lie bialgebroid in this case is 
$(T^{\pi_1}M\sdp(T^{\pi_2}M)^*,\,(T^{\pi_1}M)^*\pds T^{\pi_2}M)$ 
where the semi--direct structures are also defined by the 
connections.
\end{ex}

\begin{ex}\rm
Let $(S;H,V;M)$ be a double Lie groupoid and write $D = A_2S$. 
Consider the cotangent double groupoid $(T^*S;A_V^*S,A_H^*S;A^*K).$ 
The double Lie algebroid of $T^*S$ can be identified canonically 
with $(T^*D;D\duer A(A_V^*S),D\duer A(A_H^*S); A^*K)$, the
cotangent double of $D$ as shown in the top face of Figure
\ref{fig:cottrip}(a). 
\end{ex}

\newcommand{\noopsort}[1]{} \newcommand{\singleletter}[1]{#1} \def\cprime{$'$}
  \def\cprime{$'$}

\end{document}